\documentclass[11pt]{article}
\usepackage[a4paper,margin=0.9in]{geometry}
\usepackage{amsmath,amssymb,mathtools,mathrsfs}
\usepackage{graphicx}
\usepackage{setspace}
\usepackage[ruled,vlined,linesnumbered]{algorithm2e}
\usepackage{algpseudocode}
\usepackage{float}
\usepackage{amsthm}
\usepackage{caption}
\usepackage{subcaption}
\usepackage{booktabs,array,multirow}
\usepackage{etoolbox}
\makeatletter
\undef{\c@author}
\undef{\theauthor}
\makeatother
\usepackage[american]{babel}
\usepackage{csquotes}
\usepackage[style=apa,backend=biber]{biblatex}
\usepackage{hyperref}
\hypersetup{colorlinks=true,linkcolor=blue,citecolor=blue,urlcolor=blue}
\usepackage[capitalise,noabbrev]{cleveref}
\DeclareFieldFormat{url}{\ifhyperref{\href{#1}{\url{#1}}}{\url{#1}}\addperiod}
\DeclareFieldFormat{doi}{\ifhyperref{\href{https://doi.org/#1}{\url{https://doi.org/#1}}}{https://doi.org/#1}\addperiod}
\usepackage{chngcntr}
\counterwithout{figure}{section}
\counterwithout{table}{section}
\newtheorem{theorem}{Theorem}[section]
\newtheorem{proposition}[theorem]{Proposition}
\newtheorem{lemma}[theorem]{Lemma}
\newtheorem{assumption}{Assumption}
\newtheorem{definition}[theorem]{Definition}
\newtheorem{corollary}[theorem]{Corollary}
\newtheorem{remark}[theorem]{Remark}

\def\mb{\mathbb}

\def\mf{\mathbf}
\def\bs{\boldsymbol}
\def\E{\mb{E}}
\newcommand{\?}[1]{#1}
\newcommand{\xref}[1]{\Cref{#1}}
\newcommand{\tabref}[1]{\Cref{#1}}
\newcommand{\stmxref}[2]{\ref{#1}}
\newcommand{\xlabel}[1]{\label{#1}}
\newcommand{\xnoindent}{\noindent}
\newcommand{\xq}{}
\newcommand{\xs}{}
\newcommand{\nf}{}
\newcommand{\ubrk}{\allowbreak}
\newcommand{\brk}{\allowbreak}
\newcommand{\svert}{\vert}
\newcommand{\tmin}{-}
\newcommand{\tsup}[1]{\textsuperscript{#1}}
\newenvironment{fsrc}{}{}
\newcolumntype{P}[1]{>{\raggedright\arraybackslash}p{#1}}
\newcolumntype{L}{l}
\newcolumntype{C}{c}
\crefname{assumption}{Assumption}{Assumptions}
\Crefname{assumption}{Assumption}{Assumptions}

\title{An alternating direction method of multipliers for utility-based shortfall risk portfolio optimization}
\author{Rufeng Xiao\textsuperscript{1,$\dagger$}\qquad Zhiping Li\textsuperscript{2,$\dagger$}\qquad Rujun Jiang\textsuperscript{1,*}}
\date{}

\begin{document}
\maketitle
\begingroup
\begin{NoHyper}
\renewcommand\thefootnote{1}\footnotetext{School of Data Science, Fudan University, Shanghai 200433, P. R. China. Emails: rfxiao24@m.fudan.edu.cn; rjjiang@fudan.edu.cn}
\renewcommand\thefootnote{2}\footnotetext{Department of Industrial and Systems Engineering, University of Minnesota, Minneapolis, Minnesota 55455, United States. Email: li004950@umn.edu}
\renewcommand\thefootnote{$\dagger$}\footnotetext{Equal contributions.}
\renewcommand\thefootnote{*}\footnotetext{Corresponding author.}
\end{NoHyper}
\endgroup
\onehalfspacing

\begin{abstract}

Utility-based shortfall risk (UBSR), a convex risk measure sensitive to tail losses, has gained popularity in recent years. However, research on computational methods for UBSR optimization remains relatively scarce. In this paper, we propose a fast and scalable algorithm for the UBSR-based portfolio optimization problem. Leveraging the Sample Average Approximation (SAA) framework, we reformulate the problem as a block-separable convex program and solve it efficiently via the alternating direction method of multipliers (ADMM). In the high-dimensional setting, a key challenge arises in one of the subproblems---a projection onto a nonlinear feasibility set defined by the shortfall-risk constraint. We propose a semismooth Newton algorithm that employs an implicit-function transformation of semismooth functions to reduce the problem to a univariate equation involving the Lagrange multiplier. This algorithm achieves local superlinear convergence under mild regularity conditions and global convergence for the common loss functions considered in our experiments. Theoretical convergence guarantees of the proposed algorithms are established, and numerical experiments demonstrate a substantial speedup over state-of-the-art solvers, particularly in high-dimensional regimes.

\end{abstract}


\section{Introduction}

Accurately assessing and quantifying the risks associated with financial assets is of fundamental importance in quantitative finance and financial risk management. One of the most widely adopted risk metrics in both academic research and industry practice is \textit{Value at Risk} (VaR), which at confidence level $\alpha$ corresponds to the $(1 - \alpha)$-quantile of the loss distribution and can be interpreted as the minimum capital required to ensure that the probability of incurring a loss does not exceed $\alpha$~\parencite{duffie1997overview}. Despite its popularity, VaR ignores the magnitude of losses beyond the threshold, lacks convexity in portfolio weights, and also poses computational challenges in certain settings. To overcome these shortcomings, the axiomatic framework of coherent risk measures was introduced~\parencite{artzner1999coherent}, requiring subadditivity, translation invariance, positive homogeneity, and monotonicity. A prominent example is \textit{Conditional Value at Risk} (CVaR), defined as the expected loss conditional on exceeding the VaR at level $\alpha$, which provides a more comprehensive description of tail risk and satisfies all coherence axioms, hence yielding a convex function in portfolio weights~\parencite{rockafellar2000optimization, rockafellar2002conditional, pflug2000some, uryasev2000conditional}.

The \textit{utility-based shortfall risk} (UBSR), proposed as an alternative to VaR, has gained increasing attention in recent years and is formulated as a convex risk measure~\parencite{follmer2002convex}. In contrast to coherent risk measures, UBSR does not impose the requirements of positive homogeneity and subadditivity, while retaining convexity and economic interpretability. Given a loss function and a threshold, the UBSR of a financial position can be interpreted as the minimum amount of capital that must be added to ensure that the expected loss remains below the threshold. From this perspective, both VaR and UBSR can be viewed as capital adequacy requirements.

In particular, UBSR possesses several advantages over the conventional risk metrics, VaR and CVaR. Most notably, it offers greater flexibility in practical applications, as the choice of the loss function can reflect varying attitudes of investment managers toward financial risk. Unlike VaR, UBSR is a convex risk measure that is sensitive to the severity of tail losses. Compared to CVaR, UBSR remains invariant under randomization, as discussed in~\textcite[Remark 2.2 (4)]{dunkel2010stochastic}. Comprehensive treatments of UBSR are available; see \textcite{giesecke2008measuring} and \textcite{dunkel2010stochastic}.

While the theoretical properties of UBSR have been extensively studied, the development of computational methods for UBSR remains relatively underexplored. The existing literature~\parencite{dunkel2010stochastic, hu2018utility} primarily focuses on the estimation of UBSR by reformulating it as a root-finding problem for a nonlinear function, rather than addressing the constrained UBSR minimization problem. In particular,~\textcite{dunkel2010stochastic} employ the Robbins--Monro algorithm~\parencite{robbins1951stochastic} to solve the associated stochastic root-finding problem, while~\textcite{hu2018utility} utilize Monte Carlo-based root-finding methods within the Sample Average Approximation (SAA) framework.~\textcite{gupte2024optimization} derive the expression of the UBSR gradient under a smooth parameterization and design a stochastic projected gradient descent algorithm for UBSR estimation and constrained UBSR minimization under the Stochastic Approximation (SA) framework. Similarly,~\textcite{hegde2024online} propose a stochastic gradient descent method and extend the problem to an online setting.

Despite their computational efficiency, existing stochastic gradient methods remain challenged by high-dimensional portfolio optimization with general constraints, such as minimum expected return requirements, due to costly projection steps. Furthermore, the estimation of stochastic gradients requires strong regularity assumptions, and the algorithm is specifically designed for online learning environments, thus rendering it unsuitable for direct application in the SAA framework~\parencite{pagnoncelli2009sample, shapiro2021lectures, kim2014guide}. This SAA perspective is further motivated by practical portfolio implementation. In many applications, the problem is solved repeatedly on a rolling-window basis. This is particularly relevant when market conditions undergo structural changes, because older observations may quickly become less informative. Therefore, the computational goal is to solve the finite-sample problem induced by the current data window efficiently and accurately.

In this paper, we propose an efficient algorithm based on the Alternating Direction Method of Multipliers (ADMM) for solving the UBSR portfolio optimization problem within the SAA framework---a practical and computationally tractable setting for this class of risk measures. The portfolio selection problem~\parencite{black1992global, cui2018portfolio, chen2023distributionally} is a classical topic in quantitative finance, with its theoretical roots tracing back to Markowitz's seminal work on portfolio theory~\parencite{markowitz1952portfolio}.

We reformulate the UBSR-based SAA problem into an equivalent block-separable structure and develop an ADMM-based solution strategy. The ADMM algorithm is renowned for its convergence guarantees and scalability, and has seen wide application in fields such as machine learning and large-scale decision-making~\parencite{xiao2023unified, deng2025conic}; see also~\textcite{boyd2010distributed} and \textcite{lin2022alternating} for a comprehensive survey.

To the best of our knowledge, our algorithm is the first to be specifically designed to exploit the problem structure of UBSR-based portfolio optimization within the SAA framework, and it is accompanied by rigorous convergence guarantees. The proposed algorithm decomposes the original problem into two subproblems. Specifically, the first subproblem reduces to a convex quadratic program, which can be efficiently solved.
The second subproblem is a projection onto a convex set characterized by complex nonlinear constraints, which constitutes a major challenge in UBSR-based optimization. To address this, we develop a semismooth Newton method that solves a univariate equation via a semismooth implicit function theorem. We establish the local and global convergence of the proposed (G-)semismooth Newton methods for solving the projection subproblem under mild assumptions.
We provide an efficient implementation of the proposed algorithm, which achieves significant speedups over state-of-the-art solvers on both synthetic and real-world datasets, particularly in high-dimensional settings.

To provide a clearer perspective, \tabref{tab:ubsr_comparison} summarizes the existing computational studies on UBSR and highlights the gap our work aims to fill. As illustrated, prior research has largely concentrated on evaluating UBSR or employing stochastic gradient methods constrained by simple projection-friendly settings. In contrast, our proposed approach explicitly solves the constrained SAA optimization problem at scale.

\begin{table}[t]
\caption{\xlabel{tab:ubsr_comparison}Comparison with existing computational studies on utility-based shortfall risk (UBSR).}

\begin{tabular}{P{60pt}P{60pt}P{60pt}P{60pt}LP{60pt}}
\toprule
Reference & Problem & Method & Constraints & Framework & Key Distinction \\
\midrule
\textcite{dunkel2010stochastic}
& Fixed-position UBSR evaluation
& Robbins--Monro root-finding
& None
& SA
& Risk evaluation only \\
\hline

\textcite{hu2018utility}
& UBSR evaluation via scalar reformulation
& Monte Carlo root-finding
& General constraints
& SAA
& Risk evaluation only \\
\hline

\textcite{gupte2024optimization}
& Gradient-based UBSR optimization
& Stochastic gradient descent
& Projection-friendly only
& SA
& Smooth loss function \\
\hline

\textcite{hegde2024online}
& Online UBSR optimization
& Online stochastic gradient
& Projection-friendly only
& SA
& Online setting \\
\hline

This paper
& Constrained UBSR portfolio optimization
& ADMM + semismooth Newton
& Additional linear constraints
& SAA
& Structure-exploiting decomposition \\
\bottomrule
\end{tabular}

\end{table}

The remainder of the paper is organized as follows.  \xref{sec:prelim} introduces the theoretical preliminaries, including the semismooth implicit function theory, which forms the foundation of our algorithmic design.  \xref{sec:admm} presents the ADMM framework for UBSR-based portfolio optimization, detailing the SAA formulation and the decomposition of the problem into tractable subproblems.  \xref{sec:projection} addresses the challenging UBSR constraint projection problem, reformulating it as a constrained projection problem and developing an efficient semismooth Newton method for its solution.  \xref{sec:experiments} provides numerical results on both synthetic and real-world datasets, including UBSR projection and portfolio optimization tasks, demonstrating the effectiveness and robustness of the proposed approach.

\xnoindent\textbf{Notation.}~In this paper, boldface lowercase letters such as $\bs{x}$ denote vectors, while regular lowercase letters such as $x$ represent scalars or components of vectors. Boldface uppercase letters such as $\bs{F}$ and $\bs{G}$ denote vector-valued functions and regular uppercase letters such as $F$ and $H$ denote scalar-valued functions. For $\bs{x}\in\mb{R}^n$, $\| \bs{x}\| =\sqrt{\bs{x}^\top\bs{x}}$.
Let $\bs{0}_n$ and $\bs{1}_n$ denote the $n$-dimensional all-zero and all-one vectors, respectively, and let $\mf{I}_n$ denote the $n$-dimensional identity matrix. When the dimension is clear from context, the dimension subscript may be omitted. For a scalar $x \in \mb{R}$, denote $x_+ = \max\{x, 0\}$ and $x_- =\max\{-x,0\}$. For a positive integer $m$, we use $[m]$ to denote the set $\{1, 2, \ldots, m\}$. Let $l(\cdot): \mb{R}\rightarrow\mb{R}$ denote the univariate loss function. For a vector $\bs{x}\in\mb{R}^n$, denote $L(\bs{x}) = \sum^n_{i=1}l(x_i)$ as the sum of component-wise loss. If $l(\cdot)$ is differentiable, let $\nabla L(\bs{x})$ denote its gradient, namely, $\nabla L(\bs{x}) = [l'(x_1), l'(x_2), \ldots, l'(x_n)]^\top$, where $l'(\cdot)$ denotes the derivative of the loss function $l(\cdot)$. The operator $\mathcal{P}$ denotes the orthogonal projection; for example, $\mathcal{P}_{\mathcal{W}}(\bs{x}) = \arg\min_{\bs{u}\in\mathcal{W}}\| \bs{u}-\bs{x}\| $.

\section{Preliminaries}\xlabel{sec:prelim}

In this section, we first introduce Utility-Based Shortfall Risk (UBSR) as an alternative risk measure to Value at Risk (VaR), motivated by VaR's inability to capture tail risk and its failure to satisfy the properties of a convex risk measure. We then present the concept of G-semismoothness and the implicit function theorem for semismooth functions. These foundational concepts and assumptions serve as the theoretical basis for analyzing the ADMM algorithm applied to portfolio optimization problems in subsequent sections.

\subsection{Utility-Based Shortfall Risk}

The definition of UBSR is closely related to the loss function $l(\cdot)$ and the risk level $\lambda$, which is a prespecified parameter.
\begin{definition}[{Utility-Based Shortfall Risk~\parencite{follmer2002convex}}]
Let $X$ be a bounded integrable random variable representing a financial position. The Utility-Based Shortfall Risk with the loss function $l(\cdot)$ and the risk level $\lambda$ is defined as follows:
\begin{equation}
\xlabel{eq:UBSR def}
\text{SR}_{\lambda} (X) =\inf \left\{ t \in \mb{R}:~ \mathbb{E}\left[ l(-X - t) \right] \leq \lambda \right\}.
\end{equation}
\end{definition}

By comparing the definitions of VaR and UBSR, we observe that VaR is a special case of UBSR~\parencite{hu2018utility}.
Indeed, choosing the loss function as $l(x) = \mathbb{I}_{\{x \geq 0\}}$, where $\mathbb{I}$ denotes the indicator function, yields $\text{SR}_{\lambda}(X)=\text{VaR}_{\lambda}(X)$.
The shortcomings of VaR can be partly attributed to the properties of the loss function $l(x) = \mathbb{I}_{\{x \geq 0\}}$, such as being almost everywhere constant and non-convex within its domain. The introduction of the loss function $l(\cdot)$ fundamentally alters the properties of the risk measure. Therefore, we state the following fundamental assumption on the loss function that defines UBSR throughout this paper.
\begin{assumption}\xlabel{ass:loss monotone and convex}
The loss function $l(\cdot)$ in~\eqref{eq:UBSR def} is nondecreasing and convex, and $\inf_x l(x) < \lambda$.
\end{assumption}

We next analyze several desirable properties of UBSR, which depend closely on the characteristics of the loss function. We emphasize that \xref{ass:loss monotone and convex} aligns with both common sense and fundamental economic principles. Intuitively, the larger the negative position, the greater the loss should be. This implies that the loss function should increase monotonically. In addition, convexity captures the principle of increasing marginal loss, meaning that each additional unit of position loss results in a larger incremental loss as the size of the losing position grows. The lower bound assumption is essential to ensure the well-definedness of UBSR.
Moreover, under \xref{ass:loss monotone and convex}, UBSR is a convex risk measure, first introduced in~\textcite{follmer2002convex} as an alternative to the coherent risk measure proposed in~\textcite{artzner1999coherent}. For further properties of UBSR, readers are referred to~\textcite{giesecke2008measuring, dunkel2010stochastic}.

\subsection{Preliminaries on {\?{ADMM}}}\xlabel{sec: pre ADMM}

In this subsection, we briefly review the standard alternating direction method of multipliers (ADMM). Standard ADMM is reviewed extensively in the literature; see \textcite{boyd2010distributed} and \textcite{lin2022alternating}. Consider a convex optimization problem with linear constraints given by
\begin{equation*}
\min_{\bs{x},\bs{y}} \ f(\bs{x}) + g(\bs{y})
\quad \text{s.t.} \quad A\bs{x} + B\bs{y} = \bs{c}.
\end{equation*}
The associated augmented Lagrangian function is formulated as
\begin{equation*}
\mathcal{L}_\sigma(\bs{x},\bs{y},\bs{\lambda})
=
f(\bs{x})+g(\bs{y})
+ \bs{\lambda}^\top(A\bs{x}+B\bs{y}-\bs{c})
+ \frac{\sigma}{2}\| A\bs{x}+B\bs{y}-\bs{c}\| ^2,
\end{equation*}
where $\bs{\lambda}$ denotes the vector of dual variables, or equivalently the Lagrange multipliers, associated with the linear equality constraint, and $\sigma>0$ is the augmented Lagrangian penalty parameter. In ADMM, the choice of $\sigma$ strongly influences the conditioning of the subproblems and the balance between primal and dual residuals. Consequently, it can substantially affect practical convergence speed; see \textcite[Section 3.4.1]{boyd2010distributed}.

For the above two-block convex problem, the standard ADMM iterations are given by
\begin{equation*}
\begin{aligned}
\bs{x}^{k+1} &\coloneq  \arg\min_{\bs{x}} \mathcal{L}_\sigma(\bs{x},\bs{y}^k,\bs{\lambda}^k), \quad
\bs{y}^{k+1} \coloneq  \arg\min_{\bs{y}} \mathcal{L}_\sigma(\bs{x}^{k+1},\bs{y},\bs{\lambda}^k), \\
\bs{\lambda}^{k+1} &\coloneq  \bs{\lambda}^k + \sigma(A\bs{x}^{k+1}+B\bs{y}^{k+1}-\bs{c}).
\end{aligned}
\end{equation*}
Under standard convexity assumptions, the sequence generated by these iterations converges to a KKT point; see \textcite[Theorem 3.1]{lin2022alternating}. This framework forms the basis of our subsequent algorithmic development.

\subsection{Preliminaries on Semismooth Implicit Function Theorems}
In many scenarios, the loss function may not be a $C^\infty$ function. For example, it may be a piecewise polynomial function. In this subsection, we define semismooth functions and discuss some of their basic properties, which will be important for further analysis in the following sections. For more detailed information, the reader may refer to~\parencite{Qi1993semismooth, clarke1990nonsmooth, gowda2004inverse, li2018efficiently}.

\begin{definition}[G-semismoothness (resp. semismoothness)~{\parencite[Theorem 2.3]{Qi1993semismooth}},
{\parencite[Definition 2]{gowda2004inverse}}]
Let $\bs{F}:\Omega \to \mb{R}^m$ be a locally Lipschitz continuous function on the open set $\Omega$. We say that $\bs{F}$ is G-semismooth at $\bs{x} \in \Omega$ with respect to a nonempty, compact-valued, and upper semicontinuous set-valued mapping $\mathcal{K}$ if for any $\bs{d} \to \bs{0}$ and $J \in \mathcal{K}(\bs{x}+\bs{d})$, it holds that
\begin{equation*}
\| \bs{F}(\bs{x}+\bs{d}) - \bs{F}(\bs{x}) - J\bs{d}\|  = o(\| \bs{d}\| ) \quad \text{as } \bs{d} \to \bs{0}.
\end{equation*}
Moreover, $\bs{F}$ is strongly G-semismooth at $\bs{x} \in \Omega$ with respect to $\mathcal{K}$ if, for any $\bs{d} \to \bs{0}$ and $J \in \mathcal{K}(\bs{x}+\bs{d})$, we have
\begin{equation*}
\| \bs{F}(\bs{x}+\bs{d}) - \bs{F}(\bs{x}) - J\bs{d}\|  = O(\| \bs{d}\| ^2) \quad \text{as } \bs{d} \to \bs{0}.
\end{equation*}

Furthermore, if $\bs{F}:\Omega \to \mb{R}^m$ is both (strongly) G-semismooth at $\bs{x} \in \Omega$ with respect to $\mathcal{K}$ and directionally differentiable at $\bs{x}$, then it is (strongly) semismooth at $\bs{x} \in \Omega$ with respect to $\mathcal{K}$.
\end{definition}

\begin{remark}
In some cases, $\bs{F}$ is simply called G-semismooth (resp. semismooth) if the set-valued function $\mathcal{K}$ is replaced with its Clarke subdifferential $\partial \bs{F}$ since the definition of G-semismoothness (resp. semismoothness) requires that $\bs{F}$ is locally Lipschitz and thus $\partial \bs{F}$ is well-defined.
\end{remark}

We now present the chain rule for the composition of G-semismooth functions, affirming that such compositions remain G-semismooth. This result is a direct extension of the classic chain rule for semismooth mappings; see, e.g.,~\textcite[Theorem 7.5.17]{facchinei2003finite}.
\begin{lemma}[{\textcite[Theorem 7.5.17]{facchinei2003finite}}]
\xlabel{lem:chain G-semismooth}
Let $\bs{G}: \mathbb{R}^n \to \mathbb{R}^m$ and $\bs{F}: \mathbb{R}^m \to \mathbb{R}^k$ be G-semismooth functions with respect to $\mathcal{K}_{\bs{G}}$ and $\mathcal{K}_{\bs{F}}$, respectively.
Assume that both $\mathcal{K}_{\bs{G}}$ and $\mathcal{K}_{\bs{F}}$ are compact-valued and upper semicontinuous.
For any $\bs{x}\in\mathbb{R}^n$, it follows that
\begin{equation*}
\mathcal{K}_{\bs{H}}(\bs{x})
\coloneq  \{\, \bs{A}_{\bs{F}} \bs{A}_{\bs{G}} \mid \bs{A}_{\bs{F}} \in \mathcal{K}_{\bs{F}}(\bs{G}(\bs{x})),\; \bs{A}_{\bs{G}} \in \mathcal{K}_{\bs{G}}(\bs{x}) \,\}
\end{equation*}
is also compact-valued and upper semicontinuous.
Let $\bs{H} = \bs{F} \circ \bs{G}: \mathbb{R}^n \to \mathbb{R}^k$.
Then $\bs{H}$ is G-semismooth with respect to $\mathcal{K}_{\bs{H}}$.
Moreover, if both $\bs{F}$ and $\bs{G}$ are strongly G-semismooth, then $\bs{H}$ is strongly G-semismooth with respect to $\mathcal{K}_{\bs{H}}$.
\end{lemma}

Finally, we present the implicit function theorem for semismooth functions. This theorem is crucial for the analysis of the G-semismooth Newton method used in the projection subproblem of the ADMM algorithm introduced below.

\begin{theorem}
{(Semismooth Implicit Function Theorem~\parencite[Theorem 2, Corollary 5, and Theorem 4]{gowda2004inverse}).}\xlabel{thm:implicit function}
Let $\bs{F}:\mb{R}^n\times\mb{R}^m\to\mb{R}^m$ be semismooth with respect to its Clarke generalized Jacobian $\partial\bs{F}$ on an open neighborhood of $(\bs{x}^0,\bs{y}^0)$, and suppose that $\bs{F}(\bs{x}^0,\bs{y}^0)=\bs{0}$. Write each element of $\partial\bs{F}$ as $[A_{\bs{x}}\ A_{\bs{y}}]$ and assume that every $A_{\bs{y}}$ in $\partial\bs{F}(\bs{x}^0,\bs{y}^0)$ is nonsingular. Then there exist an open neighborhood $\mathcal{V}\subset\mb{R}^n$ of $\bs{x}^0$ and a unique local function $\bs{f}:\mathcal{V}\to\mb{R}^m$ satisfying $\bs{f}(\bs{x}^0)=\bs{y}^0$ and $\bs{F}(\bs{x},\bs{f}(\bs{x}))=\bs{0}$ for every $\bs{x}\in\mathcal{V}$. Moreover, $\bs{f}$ is G-semismooth with respect to the compact-valued and upper semicontinuous mapping
\begin{equation*}
\mathcal{K}(\bs{x})=\left\{-(A_{\bs{y}})^{-1}A_{\bs{x}}:[A_{\bs{x}}\ A_{\bs{y}}]\in\partial\bs{F}(\bs{x},\bs{f}(\bs{x}))\right\}.
\end{equation*}
\end{theorem}

\begin{remark}
Theorem 2 of \textcite{gowda2004inverse} gives the general inverse-function result, including the inverse-matrix mapping and inheritance of semismoothness. Corollary 5 specializes this result to the Clarke generalized Jacobian, for which nonsingularity replaces the orientation and index conditions. Applying this specialization to the extended function $\widetilde{\bs{F}}(\bs{x},\bs{y})=(\bs{x},\bs{F}(\bs{x},\bs{y}))$ and following the construction of Theorem 4 yields the desired result above.
\end{remark}

\section{{\?{ADMM}} for {\?{UBSR}}}\xlabel{sec:admm}

In this section, we propose an ADMM algorithm based on  \xref{sec: pre ADMM} to solve the UBSR-based portfolio optimization problem. In particular, we first reformulate the original risk minimization model using the SAA technique. This reformulation leads to a problem with a block-separable structure, which can be effectively leveraged by the ADMM framework. We then present the resulting subproblems in detail and discuss efficient solution strategies for each. This section highlights how the ADMM algorithm can be tailored to address high-dimensional, real-world optimization problems, where both subproblems can be solved efficiently due to their favorable structure.

\subsection{Problem Formulation}

In this subsection, we consider an extended version of the portfolio optimization problem under the classical Markowitz framework \ubrk \parencite{markowitz1952portfolio} by incorporating the UBSR as a risk measure. The optimization problem takes the following form:
\begin{equation}
\xlabel{p:original optimization}
\min_{\bs{w} \in \mb{R}^n} \quad (1-\alpha) ~\text{SR}_\lambda(X^\top\bs{w}) -\alpha~\mb{E}[X^\top\bs{w}] \qquad
\text{s.t.} \quad  \bs{w} \in \mathcal{W}, \mb{E}[X^\top\bs{w}]\geq R_0,
\end{equation}
where $X$ is a random vector representing the returns of $n$ financial assets, $\bs{w} \in \mb{R}^n$ denotes the portfolio weights, and $\mathcal{W}$ is a set collecting constraints of the weights. The financial position is parameterized by $\bs{w}$, such that the portfolio return is given by the linear combination $X^\top \bs{w}$. The parameter $\alpha \in [0,1)$ measures the risk aversion of the investor. Note that when $\alpha = 0$, the problem reduces to the portfolio optimization problem studied in \textcite{hu2018utility, hegde2024online}.
The proposed model constitutes a convex optimization problem that seeks to balance the expected return and the risk, measured by UBSR, while ensuring that the expected return satisfies the constraint $\mb{E}[X^\top \bs{w}] \geq R_0$, with $R_0$ denoting the predefined minimum acceptable return. A related problem has been studied in the literature~\parencite{cattaruzza2024exact}, whereas our formulation adopts UBSR as the risk measure.
We reformulate problem~\eqref{p:original optimization} in two steps. First, following the reformulation in \textcite[Section 2.4]{hu2018utility}, we introduce an auxiliary variable $t \in \mb{R}$ to decouple the UBSR evaluation from the objective. By the definition of UBSR in \eqref{eq:UBSR def}, problem~\eqref{p:original optimization} is equivalent to the following stochastic optimization problem:
\fontsize{7pt}{10.5pt}\selectfont
\begin{equation*}
\min_{\bs{w} \in \mb{R}^n,~t\in \mb{R}} (1-\alpha)~t - \alpha~\mb{E}[X^\top \bs{w}] \qquad
\text{s.t.} \quad \bs{w} \in \mathcal{W}, ~\mb{E}[X^\top\bs{w}] \geq R_0,~\mb{E} \left[ l(-X^\top \bs{w}-t) \right] \leq \lambda.
\end{equation*}\normalsize
Second, we apply the Sample Average Approximation (SAA) scheme to replace the population expectations with their empirical averages. This yields the following data-driven approximation
\fontsize{7pt}{10.5pt}\selectfont
\begin{equation}
\xlabel{p:saa optimization}
\min_{\bs{w} \in \mb{R}^n,~t\in \mb{R}} (1-\alpha)~t - \alpha~\bs{\mu}^\top \bs{w} \qquad
\text{s.t.} \quad \bs{w} \in \mathcal{W}, ~\bs{\mu}^\top\bs{w} \geq R_0,~\frac{1}{m}\sum_{i=1}^m l(-\xi_i^\top\bs{w}-t) \leq \lambda,
\end{equation}\normalsize
where $\xi_i$ are $m$ independent and identically distributed (i.i.d.) samples from $X$, $\bs{\mu}$ is the sample mean, and $t$ is introduced as in \eqref{eq:UBSR def}.
This approximation admits convergence guarantees, which can be found in~\textcite{shapiro2013consistency}.
We make the following basic assumptions on the weight set $\mathcal{W}$.
\begin{assumption}\xlabel{ass:W properties}
The set $\mathcal{W} \subseteq \mb{R}^n$ is a nonempty, closed, and convex set that is projection-friendly, meaning that for any $\bs{v} \in \mb{R}^n$, the Euclidean projection of $\bs{v}$ onto $\mathcal{W}$ can be computed efficiently. Moreover, there exists $\bar{\bs{w}}\in\operatorname{ri}(\mathcal{W})$ such that $\bs{\mu}^\top\bar{\bs{w}}>R_0$.
\end{assumption}

This assumption is typically reasonable in practice. For example, $\mathcal{W}$ can be chosen as the standard simplex $\mathcal{S} \coloneq  \bigm\{ \left.\bs{w} \in \mb{R}^n \,\right\vert \, \sum_{i=1}^n w_i = 1,\ w_i \ubrk \geq 0\ \forall i \bigm\}$, which corresponds to a portfolio weight set in a market where short selling is not allowed.
The ADMM framework presented later in Algorithm \stmxref{alg:UBSR optimization}{1} can also be applied to another classical portfolio optimization problem with the UBSR serving as a constraint. This problem involves maximizing the expected utility of a financial position subject to a prespecified upper bound on its UBSR. The reader may refer to Section~B.4 of the supplementary material for detailed implementation.

\subsection{{\?{ADMM}} for {\?{UBSR}} Optimization}

In this subsection, we present the detailed ADMM algorithmic framework and implementation steps for solving the UBSR optimization problem~\eqref{p:saa optimization}. The core challenge lies in efficiently handling the UBSR-induced constraint, which is expressed as a sample-average evaluation of a loss function applied to affine transformations of the portfolio return.

Let $\bs{R} = \left[
\begin{matrix}
\xi_1, \xi_2, \ldots, \xi_m
\end{matrix}
\right]^\top$ denote the sample return matrix of dimension $m\times n$. To address the risk minimization problem~\eqref{p:saa optimization}, we introduce an auxiliary variable $\bs{z} \in \mb{R}^m$ by setting
\begin{equation*}
\bs{z} = -\bs{R} \bs{w} - t \mf{1}_m.
\end{equation*}
This reformulation allows us to decouple the projection-friendly weight constraint $\bs{w}\in\mathcal{W}$ from the UBSR-induced constraint that is not directly amenable to projection, thereby greatly simplifying the resulting subproblems.

Furthermore, to handle the minimum expected return constraint $\bs{\mu}^\top\bs{w} \ge R_0$, we introduce a nonnegative slack variable $s$ and rewrite it as
\begin{equation*}
\bs{\mu}^\top\bs{w} - s = R_0, \quad s\geq 0.
\end{equation*}
This reformulation enables the incorporation of both the equality constraint and the slack variable into the augmented Lagrangian of the ADMM framework for subsequent updates.

Based on the preparations and transformations above, problem~\eqref{p:saa optimization} is transformed into the following problem with a separable structure
\begin{equation}
\xlabel{p:saa optimization admm}
\begin{aligned}
\min_{\bs{w}, t, \bs{z}, s} \quad&(1-\alpha)~t - \alpha~\bs{\mu}^\top \bs{w} \\
\text{s.t.} \quad& \bs{w} \in \mathcal{W}, \bs{z} \in \mathcal{Z}, s\geq 0,~\bs{z} = -\bs{R} \bs{w} - t \mf{1}_m, ~\bs{\mu}^\top\bs{w} -s= R_0,
\end{aligned}
\end{equation}
where
\begin{equation*}
\mathcal{Z} = \{\bs{z}\in\mb{R}^{m}\mid \frac{1}{m}\sum_{i=1}^m l(z_i) \leq \lambda\}.
\end{equation*}
This reformulation is equivalent to the original problem~\eqref{p:saa optimization} in the sense that problems~\eqref{p:saa optimization} and \eqref{p:saa optimization admm} have the same optimal value, and their primal optimal solutions correspond one-to-one via the $(\bs{w},t)$-component.

We now present the detailed update procedure of the ADMM algorithm. First, we derive the augmented Lagrangian function of~\eqref{p:saa optimization admm}
\fontsize{7.3pt}{10.5pt}\selectfont
\begin{align*}
L_{\sigma}(\bs{w}, t, \bs{z}, s, \bs{\nu}_1, \nu_2) = &(1-\alpha)~t - \alpha~\bs{\mu}^\top \bs{w}+ \bs{\nu}_1^\top (\bs{R}\bs{w} + t \mf{1}_m + \bs{z}) + \nu_2 (\bs{\mu}^\top \bs{w} - s - R_0) \\
&+\xs  \frac{\sigma}{2} \left [ \| \bs{R}\bs{w} + t \mf{1}_m + \bs{z}\| ^2 + (\bs{\mu}^\top \bs{w} - s - R_0)^2\right],
\end{align*}\normalsize
where $\bs{\nu}_1 \in \mb{R}^m$ and $\nu_2 \in \mb{R}$ are Lagrangian multipliers associated with the equality constraints.

We then apply ADMM to solve problem~\eqref{p:saa optimization} using the following updates:
\begin{subequations}
\xlabel{eq:optimization admm update}
\begin{align}
(\bs{w}^{k+1}, t^{k+1}) &= \arg \min_{\bs{w}, t} \left \{ L_{\sigma}(\bs{w}, t, \bs{z}^k, s^k, \bs{\nu}_1^k, \nu_2^k): \bs{w} \in \mathcal{W} \right \}, \xlabel{eq:w-t update}\\
(\bs{z}^{k+1}, s^{k+1}) &= \arg\min_{{\bs{z}}, s} \left\{ L_{\sigma}(\bs{w}^{k+1}, t^{k+1}, \bs{z}, s, \bs{\nu}_1^k, \nu_2^k) : \bs{z} \in \mathcal{Z}, s \geq 0\right \}, \xlabel{eq:z-s update} \\
\bs{\nu}_1^{k+1} &= \bs{\nu}_1^{k} + \sigma(\bs{R}\bs{w}^{k+1} + t^{k+1} \mf{1}_m + \bs{z}^{k+1}), \xlabel{eq:nu1 update} \\
\nu_2^{k+1} &= \nu_2^k + \sigma  (\bs{\mu}^\top \bs{w}^{k+1} - s^{k+1} - R_0). \xlabel{eq:nu2 update}
\end{align}
\end{subequations}
At each iteration, we alternately update one block of primal variables while holding the others fixed. Once the primal updates are complete, we perform ascent on the dual variables. Therefore, each iteration requires solving two separate subproblems over the primal variables.

Since the $(\bs{w}, t)$-update~\eqref{eq:w-t update} is a quadratic program with projection-friendly constraints, we solve it using the classical accelerated projected gradient method \parencite{nesterov2018lectures}.

For the $(\bs{z}, s)$-update~\eqref{eq:z-s update}, we need to solve the following subproblem:
\begin{align*}
(\bs{z}^{k+1}, s^{k+1}) &= \arg\min_{{\bs{z}}, ~s} \left \{\vbox to 1.5pc{}\frac{\sigma}{2} \| \bs{z} - \bs{c}_1^{k+1}\| ^2_2 \right.\\&\xq\left. +\xs \frac{\sigma}{2}\left(s-(\bs{\mu}^\top \bs{w}^{k+1}-R_0+\frac{\nu^k_2}{\sigma})\right)^2: \bs{z} \in \mathcal{Z}, s \geq 0\right \},
\end{align*}
where $\bs{c}_1^{k+1}=-\bs \nu_1^k/\sigma - \bs{R}\bs{w}^{k+1} - t^{k+1} \mf{1}_m$. Noting that this problem is separable in $\bs{z}$ and $s$, we can explicitly write down the update formulas for each variable
\begin{align}
\bs{z}^{k+1} &= \mathcal{P}_{\mathcal{Z}}(\bs{c}_1^{k+1})\coloneq  \arg\min_{\bs{z}} \left\{\| \bs{z} - \bs{c}_1^{k+1}\| ^2_2: \bs{z} \in \mathcal{Z} \right\},  \xlabel{eq:projection z} \\
s^{k+1} &= \max\{\bs{\mu}^\top \bs{w}^{k+1}-R_0+ \nu_2^k/\sigma, 0\}.\notag
\end{align}
We will propose efficient algorithms for the projection in  \xref{sec:projection}.

The complete ADMM procedure for problem~\eqref{p:saa optimization} is summarized in Algorithm \stmxref{alg:UBSR optimization}{1}. The following theorem establishes the convergence of Algorithm \stmxref{alg:UBSR optimization}{1}, showing that the iterates asymptotically satisfy the KKT conditions of the separable problem~\eqref{p:saa optimization admm} for any fixed $\sigma > 0$.
\begin{theorem}
Suppose \xref{ass:loss monotone and convex,ass:W properties} hold and the solution set of problem~\eqref{p:saa optimization} is nonempty. Let $\{(\bs{w}^k, t^k, \bs{z}^k, s^k, \bs{\nu}_1^k, \nu_2^k)\}_{k\ge0}$ be the sequence generated by Algorithm \stmxref{alg:UBSR optimization}{1}. Then, there exists a globally optimal solution $(\bs{w}^{*}, t^{*})$ for the original SAA formulation~\eqref{p:saa optimization}, along with corresponding auxiliary variables $\bs{z}^{*}, s^{*}$ and dual multipliers $\bs{\nu}_1^{*}, \nu_2^{*}$, such that the limit point
$(\bs{w}^{*}, t^{*}, \bs{z}^{*}, s^{*}, \bs{\nu}_1^{*}, \nu_2^{*})$
satisfies the KKT conditions of problem~\eqref{p:saa optimization admm}, and the following asymptotic properties hold as $k \to \infty$:
\fontsize{7pt}{10.5pt}\selectfont
\begin{equation*}
(1-\alpha)(t^k-t^{*})-\alpha\bs{\mu}^\top(\bs{w}^k-\bs{w}^{*}) \to 0,
\qquad
\bs{R}\bs{w}^k+t^k\mf{1}_m+\bs{z}^k \to \mf{0}_m,
\qquad
\bs{\mu}^\top \bs{w}^k-s^k-R_0 \to 0.
\end{equation*}\normalsize
Furthermore, the following bounds imply an $O(k^{-1/2})$ rate for the absolute objective residual and for the norm of the feasibility residual:
\fontsize{7.5pt}{10.5pt}\selectfont
\begin{equation*}
\begin{aligned}
-\| \bs{\nu}^{*}\| \sqrt{\frac{C}{\sigma(k+1)}}
\le
(1-\alpha)(t^k-t^{*})-\alpha\bs{\mu}^\top(\bs{w}^k-\bs{w}^{*})
&\le
\frac{C}{k+1}+\frac{2C}{\sqrt{k+1}}\\&\xq+\xs\| \bs{\nu}^{*}\| \sqrt{\frac{C}{\sigma(k+1)}},\\
\| \bs{R}\bs{w}^k+t^k\mf{1}_m+\bs{z}^k\| ^2
+
(\bs{\mu}^\top \bs{w}^k-s^k-R_0)^2
&\le
\frac{C}{\sigma(k+1)},
\end{aligned}
\end{equation*}\normalsize
where $(\bs{\nu}^{*})^\top=[(\bs{\nu}_1^{*})^\top,\nu_2^{*}]$ and the constant $C \coloneq  \frac{1}{\sigma}\Bigl(\| \bs{\nu}_1^0-\bs{\nu}_1^{*}\| ^2+(\nu_2^0-\nu_2^{*})^2\Bigr)
+\sigma\Bigl(\| \bs{z}^0-\bs{z}^{*}\| ^2+(s^0-s^{*})^2\Bigr)$.
\end{theorem}

\begin{proof}
Since problem~\eqref{p:saa optimization} is convex and Slater's condition (i.e.,~there exists a point $\bs{w}^0$ in the relative interior of $\mathcal{W}$ and $t^0$ such that $\bs{\mu}^\top\bs{w}^0>R_0$ and $\frac{1}{m} \sum_{i=1}^m l(-\xi_i^\top\bs{w}^0-t^0) < \lambda$) holds under \xref{ass:loss monotone and convex,ass:W properties}, the assumed nonemptiness of the solution set ensures that strong duality holds and there exists at least one saddle point $(\bs{w}^{*}, t^{*}, \bs{z}^{*}, s^{*}, \bs{\nu}_1^{*}, \nu_2^{*})$ for this problem; see \textcite[Theorem A.2]{beck2017first}. For the rest of the proof, the reader may refer to~\textcite[Theorem 3.1 and Theorem 3.2]{lin2022alternating}.
\end{proof}

\mbox{}
\begin{fsrc}
\IncMargin{12pt}
\begin{algorithm}[ht!]
	\caption{ADMM for problem~\eqref{p:saa optimization}: $\bs{w}^{k+1} = {\rm ADMM}(\bs{w}^0, \bs{z}^0, \bs{\nu}_1^0, t^0, s^0, \nu_2^0, \sigma, R_0)$}
	\label{alg:UBSR optimization}
	\DontPrintSemicolon
	\textbf{Input:} $\bs{w}^0\in\mb{R}^n$, $\bs{z}^0 \in \mb{R}^m$, $\bs{\nu}_1^0 \in \mb{R}^m$, $t^0 \in \mb{R}$, $s^0 \in \mb{R}$, $\nu_2^0 \in \mb{R}$, $\sigma > 0$, $R_0 \in \mb{R}$.\;
	\textbf{Initialize:} Set $k = 0$.\;
	\For{$k = 0,1,\ldots$}{
		Solve~\eqref{eq:w-t update} to obtain $(\bs{w}^{k+1}, t^{k+1})$, and solve \eqref{eq:z-s update} to obtain $(\bs{z}^{k+1}, s^{k+1})$.\;
		Update $(\bs{\nu}_1^{k+1},\nu_2^{k+1})$ by \eqref{eq:nu1 update} and \eqref{eq:nu2 update}.\;
	}
\end{algorithm}
\end{fsrc}

\section{\texorpdfstring{Projection onto the set $\mathcal{Z}$}{Projection onto the set Z}}\xlabel{sec:projection}
In this section, we focus on the projection problem onto the set $\mathcal{Z}$ involving the nonlinear constraint \eqref{eq:projection z}, which arises within the ADMM updates. We begin by formulating the Karush--Kuhn--Tucker (KKT) conditions~\parencite[Theorem 3.78]{beck2017first} of the projection problem. A natural approach to solving this subproblem efficiently is to treat the KKT system as a nonlinear system and apply a semismooth Newton method directly. However, such a method lacks global convergence guarantees. Therefore, we explore the possibility of leveraging the favorable properties of one-dimensional (G-)semismooth Newton methods. To this end, we introduce a G-semismooth Newton method that utilizes the implicit function theorem. This section addresses the most challenging part of the algorithmic framework.

\subsection{Projection Problem Formulation}

Building on the analysis in the previous section, we observe that a key computational component in the $\bs{z}$-update step in Algorithm \stmxref{alg:UBSR optimization}{1} involves solving the following projection problem
\begin{equation}
\xlabel{p:projection on Z}
\min_{\bs{u} \in \mb{R}^m} \quad   \frac{1}{2}\|  \bs{u} - \bs{x}\| ^2 \qquad
\text{s.t. } \quad  \bs{u} \in \mathcal{Z}=\{\bs{z}\in \mb{R}^m\mid\frac{1}{m}\sum_{i=1}^ml(z_i)\leq \lambda\}.
\end{equation}
Before presenting our analysis, we introduce a key assumption to facilitate subsequent algorithm development and convergence analysis.
\begin{assumption}\xlabel{ass:l' semismooth}
The loss function $l(\cdot)$ is differentiable and nonconstant. Its derivative $l'(\cdot)$ is semismooth with respect to $\partial l'(\cdot)$.
\end{assumption}
This assumption holds for many practical scenarios, including the exponential and piecewise polynomial functions used in our numerical experiments in  \xref{sec:experiments}.
Moreover, since the projected set $\mathcal{Z}$ is a closed sublevel set of a convex function, the following result \parencite[Theorem 6.25]{beck2017first} guarantees the uniqueness of the solution to the projection problem.
\begin{lemma}\xlabel{lem:unique solution to projection}
Let $\mathcal{C}$ be a nonempty closed convex subset of $\mb{R}^n$ and let $x\in\mb{R}^n$ be a vector. Then there exists a unique projection $x^{*} \in \mathcal{C}$.
\end{lemma}

By \xref{ass:loss monotone and convex}, the infimum of the loss function $l(\cdot)$ over $\mathbb{R}$ is strictly less than $\lambda$. Thus, the Slater condition holds for problem~\eqref{p:projection on Z}.
Then, according to \textcite[Theorem 3.78]{beck2017first}, $\bs u$ is an optimal solution of problem \eqref{p:projection on Z} if and only if there exists a multiplier $\rho$ such that the following KKT conditions hold:
\fontsize{7.5pt}{10.5pt}\selectfont
\begin{equation}
\xlabel{eq:original KKT}
u_i - x_i + \frac{\rho}{m}l'(u_i) = 0,\forall i \in [m], \quad
\rho (\frac{1}{m} \sum^m_{i=1} l(u_i) - \lambda ) = 0, ~\rho \geq 0, \quad
\frac{1}{m}\sum_{i=1}^ml(u_i) \leq \lambda.
\end{equation}\normalsize
We aim to find the solution $(\rho^{*}, \bs{u}^{*})$ to the KKT system~\eqref{eq:original KKT}. If $\bs{x} \in \mathcal{Z}$, it is trivial that $\bs u^{*}=\bs x$. Our analysis will therefore focus on the nontrivial case $\bs{x} \notin \mathcal{Z}$. We need the following helpful lemma.
\begin{lemma}\xlabel{lemma: KKT multiplier}
Suppose that $(\rho^{*}, \bs{u}^{*})$ is the solution to the KKT system~\eqref{eq:original KKT}. Then $\rho^{*} > 0$ if $\bs{x} \notin \mathcal{Z}$.
\end{lemma}
\begin{proof}
Starting with the complementary slackness condition, we observe that if $\rho^{*} = 0$, the stationarity condition of the Lagrangian implies $\bs{u}^{*} = \bs{x}$, which in turn requires that $\bs{x}$ belongs to the feasible set $\mathcal{Z}$. Hence, in the case where $\bs{x} \notin \mathcal{Z}$, it necessarily follows that $\rho^{*} > 0$.
\end{proof}
Therefore if $\bs{x} \notin \mathcal{Z}$, the corresponding KKT system is given by
\begin{equation}
\xlabel{eq:projection KKT}
u_i - x_i + \frac{\rho}{m}l'(u_i) = 0,\forall i \in [m], \quad
\frac{1}{m}\sum_{i=1}^ml(u_i) = \lambda, \quad \rho >0,
\end{equation}
where the second equality follows from the complementary slackness condition under $\rho>0$, as shown in  \xref{lemma: KKT multiplier}.

\subsection{A \?{G}-semismooth Method for Solving a Univariate Equation}\xlabel{sec: sep solve KKT}

A direct approach to~\eqref{eq:projection KKT} is to apply a semismooth Newton (SSN) method to the full high-dimensional nonlinear system. Although this approach enjoys local superlinear convergence, globalization techniques, such as line-search, are typically necessary to avoid divergence from random initialization; see Section~3 of \textcite{pang1990newton}, Section~3 of \textcite{jiang1998global}, and Section~2 of \textcite{qi1999survey}. Section~A of the supplementary material presents this direct SSN approach. This motivates the question of how to retain the efficiency of semismooth Newton methods while establishing global convergence. We observe that, in contrast to its higher-dimensional counterparts, the univariate (G-)semismooth Newton method for solving scalar nonlinear equations possesses the rigorous property of global convergence under specific regularity conditions.
Theorem 2.1 in~\textcite{aravkin2019level} states that, for a univariate function that is decreasing, convex, and satisfies certain regularity conditions, the (G-)semismooth Newton method achieves at least $Q$-superlinear convergence to its root. Therefore, a promising approach with stronger convergence guarantees is to construct a transformation that reduces the high-dimensional nonlinear system into a sequence of univariate nonlinear equations, each solved via the (G-)semismooth Newton method. Recall that the simplified KKT system~\eqref{eq:projection KKT} consists of $m+1$ equations in $m+1$ unknowns. As defined in the Notation part, $L(\bs{u}) = \sum^m_{i=1} l(u_i)$ represents the aggregate loss function over the components of $\bs{u}\in\mb{R}^m$. For ease of exposition, we define the two functions
\begin{equation}
\xlabel{eq:H and G}
\begin{aligned}
&H(\rho): = L(\bs{u}(\rho))-m\lambda = \sum^m_{i=1} l(u_i(\rho)) - m\lambda, \\
&\bs{G}(\bs{u}, \rho): = \bs{u} - \bs{x} + \frac{\rho}{m}\nabla L(\bs{u}).
\end{aligned}
\end{equation}

To efficiently compute $(\rho^{*}, \bs{u}^{*})$ satisfying $H(\rho^{*}) = 0$ and $\bs{G}(\bs{u}^{*}, \rho^{*}) = \bs{0}_{m}$, we reformulate the problem as a one-dimensional root-finding task exploiting the implicit relation $\bs{G}(\bs{u}, \rho) = \bs{0}_{m}$. Under mild regularity conditions, for a given $\rho$, both $H(\rho)$ and its derivative can be evaluated through the implicit function $\bs{u}(\rho)$ and its derivative with respect to $\rho$. When $l'$ is convex, the semismooth Newton method for the $\bs{G}$-subproblem is globally convergent. For the outer equation $H(\rho)=0$, we establish local superlinear convergence in general and global convergence for common loss functions under the additional conditions given below.

In the following, we first propose a semismooth Newton method for solving $\bs{G}(\bs{u}, \rho) = \bs{0}_{m}$ and analyze its convergence. Subsequently, we will derive the properties of the univariate function $H(\rho)$ and propose an associated G-semismooth Newton method. Furthermore, we will verify that $H(\rho)$ is convex under some common loss functions, which allows us to establish the global convergence of the proposed G-semismooth Newton method.

\subsubsection{\texorpdfstring{Semismooth \?{Newton} Method for $\bs{G}$-subproblem}{Semismooth Newton Method for G-subproblem}}

In this part, we focus on establishing the convergence of the semismooth Newton method applied to the equation $\bs{G}(\bs{u}, \rho) = \bs{0}_m$ for a fixed $\rho > 0$.
Given $\bs{x}$, by Corollary 2.4 in~\textcite{Qi1993semismooth}, $\nabla L$ is semismooth under \xref{ass:l' semismooth}. Together with the fact that $\bs{G}(\bs{u}, \rho)$ is affine in $\rho$, we conclude that the function $\bs{G}(\bs{u}, \rho)$ is semismooth. We define $\bs{G}_\rho(\bs{u}) \coloneq  \bs{G}(\bs{u}, \rho)$ for a fixed $\rho$, which preserves the semismoothness property. The Clarke subdifferential of $\bs{G}_\rho(\bs{u})$ can be characterized as
\begin{equation}
\xlabel{eq:partial G_k}
\partial \bs{G}_\rho(\bs{u}) = \left\{\left. \bs{I}_m + \frac{\rho}{m} \Lambda \,\right\vert \, \Lambda \in \partial \nabla L(\bs{u}) \right\},
\end{equation}
where we recall that \brk $
\partial \nabla L(\bs{u}) = \left\{\left. \mathrm{diag}(\eta_1, \eta_2, \ldots, \eta_m) \,\right\vert \, \eta_i \in \partial l'(u_i) \text{ for all } i \right\}$.
\begin{lemma}\xlabel{lem:G-rho-nonsingular}
Suppose \xref{ass:loss monotone and convex,ass:l' semismooth} hold. Then, for any $\bs{u}\in\mathbb{R}^m$, any $\rho>0$, and any $\Lambda\in\partial\nabla L(\bs{u})$, the matrix $\Lambda$ is diagonal with nonnegative diagonal elements, and $\bs{I}_m+\frac{\rho}{m}\Lambda$ is nonsingular.
\end{lemma}
\begin{proof}
Since $l$ is a convex differentiable function, $l'$ is a nondecreasing function. By \xref{ass:l' semismooth}, $l'$ is locally Lipschitz continuous due to its semismoothness. Then for any $\Lambda \in \partial\nabla L(\bs{u})$, by Proposition 2.3 in~\textcite{jiang1995local}, we know that $\Lambda$ is a diagonal matrix with nonnegative diagonal elements, which means that $\bs{I}_m + \frac{\rho}{m}\Lambda$ is nonsingular for $\rho>0$.
\end{proof}
Thus, we are ready to apply the semismooth Newton method to solve the equation $\bs{G}_\rho(\bs{u}) = \bs{0}_m$. The detailed implementation of this procedure is summarized in Algorithm \stmxref{alg:G subproblem}{2}.

\begin{fsrc}
\IncMargin{12pt}
\begin{algorithm}[ht!]
	\caption{Semismooth Newton Method for $\bs{G}_\rho(\bs{u}) = \bs{0}_m$}
	\label{alg:G subproblem}
	\DontPrintSemicolon
	\textbf{Input:} Given $\bs{x} \in \mb{R}^m$, $\rho > 0$, initial point $\bs{u}^0 \in \mb{R}^m$.\;
	\For{$t = 0,1,\ldots$}{
		Compute $\bs{G}_\rho(\bs{u}^t) = \bs{u}^t - \bs{x} + \frac{\rho}{m} \nabla L(\bs{u}^t)$ and select $P_t \in \partial \bs{G}_\rho(\bs{u}^t)$ defined in~\eqref{eq:partial G_k}.\;
		Update $\bs{u}^{t+1} = \bs{u}^t - P_t^{-1} \bs{G}_\rho(\bs{u}^t)$.\;
	}
\end{algorithm}
\end{fsrc}

We now present the local convergence theorem for Algorithm \stmxref{alg:G subproblem}{2}, along with a complete proof to ensure self-contained exposition.
\begin{theorem}
Suppose \xref{ass:loss monotone and convex,ass:l' semismooth} hold.
Let $\bs{x}\in \mathbb{R}^m$, $\rho>0$ and let $\bs u^{*}$ be a root of
$\bs{G}_\rho(\bs{u}^{*}) = 0$.
Then there exists $\delta>0$ such that, for any initial point
$\bs{u}^0 \in N(\bs{u}^{*},\delta)$, the sequence
$\{\bs{u}^t\}_{t=0}^\infty $ generated by
Algorithm \stmxref{alg:G subproblem}{2} converges to $\bs{u}^{*}$  at a $Q$-superlinear rate.
\end{theorem}
\begin{proof}
By  \xref{lem:G-rho-nonsingular}, every $P \in \partial\bs{G}_\rho(\bs{u}^{*})$ defined in~\eqref{eq:partial G_k} is nonsingular. As a result, by Theorem 3.2 in~\textcite{Qi1993semismooth}, we obtain the desired result.
\end{proof}

Moreover, since $\bs{G}_\rho(\bs{u})$ is separable, the equation $\bs{G}_\rho(\bs{u}) = \bs{0}$ decouples into $m$ independent one-dimensional root-finding problems. When $l'$ is convex, we see that $\bs{G}_\rho(\bs{u})$ is convex. Then we use the additional convexity condition to show that Algorithm \stmxref{alg:G subproblem}{2} either converges globally or terminates in finitely many steps.
This result covers a large class of loss functions, including the exponential loss $l(x) = \exp(\beta x)$ with $\beta > 0$ and the piecewise polynomial loss $l(x) = \tfrac{1}{\eta}(x_+)^\eta$ with $\eta \geq 2$, both of which are used in our numerical experiments.

\begin{theorem}\xlabel{thm:G-converges}
Suppose \xref{ass:loss monotone and convex,ass:l' semismooth} hold,
and assume that $l'$ is convex.
Let $\bs{x}\in \mathbb{R}^m$, $\rho>0$.
Then, for any initial point $\bs{u}^0\in \mathbb{R}^m$,
the sequence $\{\bs{u}^t\}_{t=0}^\infty$ generated by Algorithm \stmxref{alg:G subproblem}{2}
either converges globally to a point $\bs{u}^{*}$ satisfying $\bs{G}_\rho(\bs{u}^{*}) = 0$
at a $Q$-superlinear rate, or terminates in finitely many steps.
\end{theorem}
\begin{proof}
Since $\bs{G}_\rho(\bs{u})$ admits the structure $\bs{G}_\rho(\bs{u}) = [\bs{G}_{\rho,1}(u_1),\ \dots,\ \ubrk \bs{G}_{\rho,m}(u_m)]^\top$, each component of the equation $\bs{G}_\rho(\bs{u}) = \bs{0}$ reduces to a one-dimensional root-finding problem. We start by establishing the existence and uniqueness of the root for each coordinate.

By  \xref{ass:loss monotone and convex,ass:l' semismooth}, the function $l'(\cdot)$ is nonnegative and nondecreasing. Consequently, for any $i \in [m]$, the function $
\bs{G}_{\rho,i}(u_i) = u_i - x_i + \frac{\rho}{m}l'(u_i)$
satisfies $\lim_{u_i \to \infty} \bs{G}_{\rho,i}(u_i) = \infty$ and $\lim_{u_i \to -\infty} \bs{G}_{\rho,i}(u_i) = -\infty$. Furthermore, $\bs{G}_{\rho,i}(u_i)$ is strictly increasing, and thus admits a unique root $u_i^{*}$ such that $\bs{G}_{\rho,i}(u_i^{*}) = 0$.

Define the index set
$
I(\bs{x}) \coloneq  \left\{ i \mid l'(x_i) = 0 \right\}
$ and $b \coloneq  \max \bigm\{ x_i \mid i \in \ubrk  [m]\setminus I(\bs{x})\bigm\}$.
Furthermore, let
\begin{equation}
\xlabel{eq: l'=0 lower}
a \coloneq  \left\{
\begin{aligned}
&-\infty && \text{if} \quad l'(u)>0 \ \forall u,\\
&\sup \left\{ u \mid l'(u) = 0 \right\}  &&\text{otherwise}.
\end{aligned}
\right.
\end{equation}
Then for any $ i \in [m]\setminus I(\bs{x}) $, the root $ u_i^{*} $ lies in the interval $ (a, b) $.
In addition, we show that each component of $\bs{G}_\rho$ satisfies the non-degeneracy condition, i.e.,~$
\sup\{g \mid g \in \partial \bs{G}_{\rho,i}(u_i^{*})\} > 0,
$
for all $i \in [m]\setminus I(\bs{x})$. Indeed, we have
$
\partial \bs{G}_{\rho,i}(u_i^{*}) = 1 + \frac{\rho}{m} \partial l'(u_i^{*}).
$
By  \xref{lem:G-rho-nonsingular}, every element of $\partial l'(u_i^{*})$ is nonnegative. Consequently, $\partial \bs{G}_{\rho,i}(u_i^{*}) \geq 1$. Hence, the non-degeneracy condition is satisfied for all $i \in [m]\setminus I(\bs{x})$. Moreover, since $l'$ is a convex function, we have that $\bs{G}_\rho$ is convex.
Since $\bs{G}_{\rho,i}$ is increasing, we cannot directly apply Theorem 2.1 in~\textcite{aravkin2019level}. However, the desired result can still be obtained through a similar proof.
\end{proof}

\subsubsection{\texorpdfstring{Properties of $H(\rho)$ and Convergence Results}{Properties of H(rho) and Convergence Results}}

In this part, we first derive a differential property of $\bs{u}$ with respect to $\rho$ based on the implicit relation defined by $\bs{G}(\bs{u}, \rho) = \bs{0}$. This analysis enables us to characterize the behavior of $\bs{u}(\rho)$ as a function of $\rho$, which is essential for establishing the semismoothness of the function $H(\rho)$ defined in \eqref{eq:H and G}. In particular, we aim to show that $H(\rho)$ inherits the semismoothness from $\bs{G}(\bs{u}, \rho)$ through the implicit function $\bs{u}(\rho)$, thereby justifying the use of a G-semismooth Newton method to solve $H(\rho) = 0$.

\begin{corollary}\xlabel{corollary: u G semismooth}
Under \xref{ass:loss monotone and convex,ass:l' semismooth}, there exists a unique G-semismooth function $\bs{u} = \bs{u}(\rho)$ with respect to $\mathcal{K}_{\bs{u}}(\rho)$ defined below satisfying $\bs{G}(\bs{u}(\rho), \rho) = \bs{0}_m$ for $\rho > 0$,
\begin{equation*}
\mathcal{K}_{\bs{u}}(\rho) = \{ -(m\bs{I}_m + \rho~\Lambda)^{-1} \nabla L(\bs{u}) \mid \Lambda \in \partial \nabla L(\bs{u}) \}.
\end{equation*}
Then, $H(\rho)$ is a G-semismooth function with respect to a set-valued mapping $\mathcal{K}_{H}(\rho)$ and
\begin{equation}
\xlabel{eq:partial H_rho}
\mathcal{K}_{H}(\rho) = \{ -\nabla L(\bs{u})^\top(m\bs{I}_m + \rho\Lambda)^{-1} \nabla L(\bs{u}) \mid \Lambda \in \partial \nabla L(\bs{u}) \}.
\end{equation}
Moreover, if the derivative of the loss function $l'(\cdot)$ is strongly semismooth, $H(\rho)$ is also strongly G-semismooth.
\end{corollary}
\begin{proof}
First, we show that the implicit function $\bs{u}(\rho)$ is G-semismooth with respect to $\mathcal{K}_{\bs{u}}(\rho)$.
We obtain that there is a unique point $\bs{u}^0$ given $\rho_0>0$ which satisfies $\bs{G}(\bs{u}^0, \rho_0)=0$ from the second paragraph in the proof of  \xref{thm:G-converges}.
Since $l'$ is a semismooth function under \xref{ass:l' semismooth}, we can conclude that $\bs{G}(\bs{u}, \rho)$ is semismooth with respect to
\begin{equation*}
\partial \bs{G}(\bs{u},\rho) = \{ \left(\bs{I}_m + \frac{\rho}{m}\Lambda ,\frac{1}{m}\nabla L(\bs{u})\right) \mid \Lambda \in \partial \nabla L(\bs{u})\}.
\end{equation*}

Let $\bs{G}_{\rho_0}(\bs{u})=\bs{G}(\bs{u}, \rho_0)$. Then, $\partial \bs{G}_{\rho_0}(\bs{u}) = \left\{ \bs{I}_m + \frac{\rho_0}{m}\Lambda \mid \Lambda \in \partial \nabla L(\bs{u}) \right\}$. By  \xref{lem:G-rho-nonsingular}, every element in $\partial \bs{G}_{\rho_0}(\bs{u}^0)$ is nonsingular.
By  \xref{thm:implicit function}, there exists an open set $\mathcal{V} \subset \mb{R}$ containing $\rho_0$ and the function $\bs{u}(\rho): \mathcal{V} \rightarrow \mb{R}^m$ is G-semismooth with respect to the set-valued function $\mathcal{K}_{\bs{u}}$ satisfying $\bs{u}(\rho_0)=\bs{u}^0$ and $\bs{G}(\bs{u}(\rho),\rho)=\bs{0}$ for every $\rho \in \mathcal{V}$. The set-valued mapping $\mathcal{K}_{\bs{u}}$ is defined as follows:
\begin{equation*}
\mathcal{K}_{\bs{u}}: \rho \rightarrow  \{ -(m\bs{I}_m + \rho\Lambda)^{-1} \nabla L(\bs{u}) \mid \Lambda \in \partial \nabla L(\bs{u}) \}.
\end{equation*}
By~ \xref{lem:chain G-semismooth}, we obtain that $H(\rho)$ is G-semismooth with respect to $\mathcal{K}_{H}(\rho)$ defined in~\eqref{eq:partial H_rho}.

Moreover, if $l'(\cdot)$ is strongly semismooth, then the proof of \textcite[Theorem 2(g)]{gowda2004inverse} can be adapted by strengthening the remainder estimate from $o(\| \cdot\| )$ to $O(\| \cdot\| ^2)$, thereby establishing the strong G-semismoothness of $\bs{u}$. It then follows from  \xref{lem:chain G-semismooth} that $H$ is also strongly G-semismooth.
\end{proof}

By  \xref{corollary: u G semismooth}, we can adopt the G-semismooth Newton method to solve $H(\rho) = 0$, as $H(\rho)$ is G-semismooth with respect to $\mathcal{K}_{H}(\rho)$. At each iteration $j$, given an input $\rho^j$, we first solve the nonlinear system $\bs{G}(\bs{u}, \rho^j) = \bs{0}_m$ to obtain the corresponding $\bs{u}^j = \bs{u}(\rho^j)$. Then we update $\rho^{j+1}$ using an element from the set $\mathcal{K}_{H}(\rho^j)$. We summarize the entire procedure in Algorithm \stmxref{alg:semismooth newton}{3}.

\begin{fsrc}
\IncMargin{12pt}
\begin{algorithm}[ht!]
	\caption{One-Dimensional G-Semismooth Newton Method for Problem~\eqref{p:projection on Z}}
	\label{alg:semismooth newton}
	\DontPrintSemicolon
	\textbf{Input:} Given $\bs{x} \not\in \mathcal{Z}$, choose initial value $\rho^0 >0$.\;
	\For{$k = 0,1,\ldots$}{
		Update $\bs{u}^k$ by solving $\bs{G}_{\rho^k}(\bs{u}) = \bs{0}_m$ using Algorithm~\ref{alg:G subproblem}.\;
		Compute $H(\rho^k) = L(\bs{u}^k) - m\lambda$ and select $h_k \in \mathcal{K}_{H}(\rho^k)$ defined in~\eqref{eq:partial H_rho}.\label{line:jacob}\;
		Update $\rho^{k+1} = \rho^k - h_k^{-1} H(\rho^k)$.\;
	}
\end{algorithm}
\end{fsrc}

Next, we provide the convergence analysis of Algorithm \stmxref{alg:semismooth newton}{3}. Before proceeding, we first present a useful lemma that establishes that every element in $\mathcal{K}_{H}(\rho)$ at the solution is nonsingular. This property plays a critical role in guaranteeing the local superlinear convergence of the G-semismooth Newton method.

\begin{lemma}\xlabel{lem:H well-defined}
Suppose \xref{ass:loss monotone and convex,ass:l' semismooth} hold and $\bs{x} \notin \mathcal{Z}$  in problem~\eqref{p:projection on Z}. Then, $H(\rho)$ is nonincreasing in $\rho \in (0, \infty)$. Furthermore, the optimal Lagrange multiplier in the KKT system~\eqref{eq:projection KKT} satisfies $\rho^{*} > 0$, with $h < 0$ for all $h \in \mathcal{K}_{H}(\rho^{*})$, where $\mathcal{K}_{H}(\rho^{*})$ is defined in \eqref{eq:partial H_rho}.
\end{lemma}

\begin{proof}
We first show that $H(\rho)$ is nonincreasing in $\rho \in (0, +\infty)$. Recall that
\begin{equation*}
H(\rho) = \sum_{i=1}^m l\bigl(u_i(\rho)\bigr) - m\lambda.
\end{equation*}
By \xref{ass:loss monotone and convex}, $l$ is nondecreasing, so the result follows if $u_i(\rho)$ is nonincreasing in $\rho\in (0, \infty)$ for all $i$. From the second paragraph in the proof of  \xref{thm:G-converges}, we know that $u_i(\rho)$ is defined as the unique solution of
\begin{equation*}
u_i - x_i + \tfrac{\rho}{m} l'(u_i) = 0.
\end{equation*}
We prove that $u_i(\rho)$ is nonincreasing in $\rho \in (0, \infty)$ for all $i$ by contradiction. Without loss of generality, assume that $\rho_2 > \rho_1 > 0$. By \xref{ass:loss monotone and convex}, $l$ is convex, so it follows that $l'$ is nondecreasing, as stated in Theorem 2.13 of \textcite{rockafellar2009variational}. If $u_i(\rho_2) > u_i(\rho_1)$, then we have $l'(u_i(\rho_2)) \geq l'(u_i(\rho_1))$, which implies that
\begin{equation*}
\tfrac{\rho_2}{m} l'(u_i(\rho_2)) \geq \tfrac{\rho_1}{m} l'(u_i(\rho_1)).
\end{equation*}
This leads to the inequality
\begin{equation*}
x_i = u_i(\rho_2) + \tfrac{\rho_2}{m} l'(u_i(\rho_2)) > u_i(\rho_1) + \tfrac{\rho_1}{m} l'(u_i(\rho_1))=x_i,
\end{equation*}
which causes a contradiction. So $u_i(\rho_2) \leq u_i(\rho_1)$, which means that $u_i(\rho)$ is nonincreasing in $\rho \in (0, \infty)$.

Next we prove that the optimal Lagrange multiplier in the KKT system~\eqref{eq:projection KKT} satisfies $\rho^{*} > 0$ with $h < 0$ for all $h \in \mathcal{K}_{H}(\rho^{*})$.
If $\bs{x} \notin \mathcal{Z}$, we know that $\rho^{*} > 0$ by  \xref{lemma: KKT multiplier}.
With $\rho^{*} > 0$ and $\bs{x} \notin \mathcal{Z}$, we claim that there exists at least one dimension $j$ of $\bs{u}^{*}$ such that $\frac{\rho^{*}}{m}l'(u^{*}_j) \ne 0$. Indeed, if this is not the case, we have $\frac{\rho^{*}}{m}l'(u^{*}_j) = 0$ for all $j$, implying that $\bs{u}^{*}$ minimizes $L(\bs u)$. However, this, together with \xref{ass:loss monotone and convex}, implies that $L(\bs{u}^{*})<m\lambda$, contradicting $\bs{x} \notin \mathcal{Z}$. Therefore, we obtain that for any $h \in \mathcal{K}_{H}(\rho^{*})$,
\begin{equation*}
h = -\sum^m_{i=1}\frac{l'(u_i^{*})^2}{m + \rho^{*}\eta_i^{*}} \leq -\frac{l'(u_j^{*})^2}{m + \rho^{*}\eta_j^{*}} < 0,
\end{equation*}
where $\text{diag}( \eta_1^{*}, \eta_2^{*}, \ldots, \eta^{*}_m) \in \partial \nabla L(\bs{u}^{*})$.
\end{proof}

This ensures that the subsequent Newton step is well-defined near $\rho^{*}$. Then we present the convergence theorem for Algorithm \stmxref{alg:semismooth newton}{3}.
\begin{theorem}
Suppose \xref{ass:loss monotone and convex,ass:l' semismooth} hold.
Let $\bs{x}\in \mathbb{R}^m$.
Then there exists $\delta>0$ such that, for any initial point
$\rho^{0}\in N(\rho^{*},\delta)$, where $\rho^{*}$ satisfies $H(\rho^{*})=0$,
the sequence $\{\rho^{k}\}_{k=0}^{\infty}$ generated by Algorithm \stmxref{alg:semismooth newton}{3} converges to $\rho^{*}$ at a $Q$-superlinear rate.
Moreover, if $H(\rho)$ is strongly G-semismooth, the convergence of $\rho^{k}$ to $\rho^{*}$ is quadratic.
\end{theorem}
\begin{proof}
By  \xref{lem:chain G-semismooth} and  \xref{corollary: u G semismooth}, we have that the set-valued function $\mathcal{K}_{H}(\rho)$ is a Newton map which is defined in \textcite[Definition 2]{klatte2018approximations}.
By  \xref{lem:H well-defined}, there exists $ M > 0 $ such that $ \svert 1/h\svert  \leq M $ for all $ h \in \mathcal{K}_{H}(\rho^{*}) $. Therefore, $ \mathcal{K}_{H} $ is weak-Newton regular at $ \rho^{*} $ \parencite[Definition 3]{klatte2018approximations}. For the rest of the proof, the reader may refer to Theorem 4 in~\textcite{klatte2018approximations}.
\end{proof}

In general, $H(\rho)$ may not be convex, so we cannot directly apply the global convergence result from Theorem 2.1 in~\textcite{aravkin2019level}. However, since $H(\rho)$ is a univariate and nonincreasing function, it naturally suggests using the bisection method for root finding to provide a global convergence guarantee. The complete algorithm is presented as Algorithm \stmxref{alg:bisec projection on Z}{4}.
We next give its convergence result.

\begin{theorem}
Suppose \xref{ass:loss monotone and convex,ass:l' semismooth} hold.
Let $\rho_U>\rho_L$ and $\epsilon>0$.
Then Algorithm \stmxref{alg:bisec projection on Z}{4} takes at most $k$ iterations to generate a point $\rho^k$ satisfying $\svert \rho^k-\rho^{*}\svert  < \epsilon$,
where $\rho^{*}$ satisfies $H(\rho^{*})=0$ and $k = \lceil\log_2\!\left(\frac{(\rho_U-\rho_L)(\rho^{*}+1)}{\epsilon}\right)\rceil$.
\end{theorem}
\begin{proof}
From  \xref{lem:unique solution to projection,lem:H well-defined}, we know that $H(\rho)$ is monotone and admits a unique root satisfying $H(\rho) = 0$. By \eqref{eq:projection KKT}, we obtain that $\rho_L = 0$ serves as a lower bound of $\rho$. Then the iteration number $k$ consists of two parts:
finding an interval $(\rho_L, \rho_U]$ that contains the root, and the binary search for an approximate root.
For the first part, since the value of $\rho_U$ is multiplied by 2 at each step, we obtain that $\log_2(\rho^{*}+1)$ serves as an upper bound. An upper bound for the second part is $\log_2(\frac{\rho_U-\rho_L}{\epsilon})$. Combining the two upper bounds we have $k = \lceil \log_2\!\left(\frac{(\rho_U-\rho_L)(\rho^{*}+1)}{\epsilon}\right)\rceil $ as desired.
\end{proof}

\mbox{}
\begin{fsrc}
\IncMargin{12pt}
\begin{algorithm}[ht!]
	\caption{Bisection Method for Problem~\eqref{p:projection on Z}}
	\label{alg:bisec projection on Z}
	\DontPrintSemicolon
	\textbf{Input:} $\bs{x} \in \mb{R}^m$ with $\bs{x} \notin \mathcal{Z}$, lower bound $\rho_L = 0$, upper bound $\rho_U>0$, and tolerance $\epsilon>0$.\;
	Set $\bs{u}^0$ by solving $\bs{G}_{\rho_U}(\bs{u}) = \bs{0}_m$ using Algorithm~\ref{alg:G subproblem}.\;
	\While{$H(\rho_U)=\sum_{i=1}^m l(u^0_i(\rho_U)) - m\lambda > 0$}{
		Update $\rho_L \gets \rho_U$ and $\rho_U \gets 2  \rho_U$.\;
		Update $\bs{u}^0$ by solving $\bs{G}_{\rho_U}(\bs{u}) = \bs{0}_m$ using Algorithm~\ref{alg:G subproblem}.\;
	}
	\For{$k = 0,1,\ldots$}{
		Set $\rho^k = (\rho_L + \rho_U)/2$.\;
		Update $\bs{u}^k$ by solving $\bs{G}_{\rho^k}(\bs{u}) = \bs{0}_m$ using Algorithm~\ref{alg:G subproblem}.\;
		Compute $H(\rho^k) = \sum_{i=1}^m l(u^k_i(\rho^k)) - m\lambda$.\;
		\If{$H(\rho^k) > 0$}{
			$\rho_L \gets \rho^k$.\;
		}
		\Else{
			$\rho_U \gets \rho^k$.\;
		}
		\If{$|\rho_U - \rho_L| < \epsilon$ \textbf{or} $H(\rho^k) = 0$}{
			\Return $\bs{u}^k$.\;
		}
	}
\end{algorithm}
\end{fsrc}

In practice, we observed that the G-semismooth Newton method consistently finds the root faster than the bisection method, as shown in the experimental section.

\subsubsection{Global Convergence for Two Common Loss Functions}

In this part, we show that $H(\rho)$ is convex under certain conditions, which further guarantees the global convergence of Algorithm \stmxref{alg:semismooth newton}{3}.

The following corollary shows that when the loss function is twice continuously differentiable and possesses specific properties, $H(\rho)$ is a convex and differentiable function.
In this case, Algorithm \stmxref{alg:semismooth newton}{3} reduces to the classical Newton's method, enjoying global $Q$-superlinear convergence.

\begin{corollary}
Suppose that \xref{ass:loss monotone and convex,ass:l' semismooth} hold, the loss function $l(\cdot)$ is twice continuously differentiable and its second derivative $l''(\cdot)$ is locally Lipschitz on $(a, +\infty)$ where $a$ is defined in \eqref{eq: l'=0 lower}. Suppose further that the loss function satisfies the following condition
\begin{equation}
\xlabel{eq:pos_cov}
3l''(u)^2 - \xi l'(u) \geq 0,
\end{equation}
for any $\xi \in \partial l''(u)$.
Then the univariate function $H(\rho)$ is convex on $\rho > 0$ under this setting. Furthermore, $H(\rho)$ is differentiable, and the term $\mathcal{K}_H(\rho^k)$ in Line~4 of Algorithm \stmxref{alg:semismooth newton}{3} can be replaced
with $\{H'(\rho^k)\}$, implying that Algorithm \stmxref{alg:semismooth newton}{3} reduces to the classical Newton's method.
The sequence $\{\rho^k\}_{k=0}^\infty$ generated by Algorithm \stmxref{alg:semismooth newton}{3} converges either finitely or $Q$-superlinearly to $\rho^{*}$ satisfying $H(\rho^{*}) = 0$, for any initial point $\rho^0 > 0$.
\end{corollary}
\begin{proof}
Since the loss function $l(\cdot)$ is twice continuously differentiable, the classical implicit function theorem \parencite[Theorem A.2]{nocedal2006numerical} implies that the derivatives of $u_i(\rho)$ and $H(\rho)$ are given by
\begin{equation*}
u_i'(\rho) = -\frac{l'(u_i)}{m + \rho l''(u_i)}, \quad \forall i \in [m], \quad \text{and} \quad
H'(\rho) = -\sum_{i=1}^m \frac{l'(u_i)^2}{m + \rho l''(u_i)} \leq 0.
\end{equation*}
Since $m + \rho l''(u_i) \geq m>0$, $u_i'(\rho)$ is always finite, implying that $u_i(\rho)$ is locally Lipschitz. Moreover, since $l''$ is locally Lipschitz, we obtain that $l''(u_i(\rho))$ is also locally Lipschitz and its Clarke subdifferential exists at every $\rho > 0$ \parencite[Proposition 2.1.2]{clarke1990nonsmooth}. By the chain rule of Clarke differential \parencite[Theorem 2.3.9]{clarke1990nonsmooth}, we can obtain that
\begin{align*}
\partial H'(\rho) &\subset  \left\{\sum_{i=1}^m \frac{l'(u_i)^2}{\left(m + \rho l''(u_i)\right)^3} \left( 3~m l''(u_i) \right.\right.\\&\xq+\xs \rho \left.\left.\left( 3l''(u_i)^2 - \xi_il'(u_i) \right) \right)\mid \xi_i \in \partial l''(u_i),~u_i = u_i(\rho)\vbox to 1.2pc{}\right\},
\end{align*}
where $u_i = u_i(\rho)$. By \eqref{eq:pos_cov} and the nonnegativity of $l''(u_i)$ under \xref{ass:loss monotone and convex}, every element in $\partial H'(\rho)$ is nonnegative. Therefore, by Proposition 2.3 in~\textcite{jiang1995local}, $H'(\rho)$ is nondecreasing with respect to $\rho$. Then by Theorem 2.13 in \textcite{rockafellar2009variational}, we know that $H(\rho)$ is convex on $(0,+\infty)$.
Similar to the proof of \xref{lem:H well-defined}, $H$ satisfies the nondegeneracy condition because $H'(\rho^{*})<0$. Recall that $H'(\rho) \leq 0$ for all $\rho >0$; thus, $H$ is nonincreasing in $\rho \in (0,\infty)$. The rest of the proof follows similarly to that of Theorem 2.1 in~\textcite{aravkin2019level}.
\end{proof}
It is straightforward to verify that the exponential loss $l(x) = \exp(\beta x)$ with $\beta > 0$ and the piecewise polynomial function $l(x) = \eta^{-1}(x_+)^{\eta}$ with $\eta \geq 3$ satisfy \eqref{eq:pos_cov}.

The preceding corollary provides a sufficient condition based on third-order information of the loss function to verify the convexity of $H(\rho)$, thereby ensuring the global convergence of the Newton sequence.
However, some prevalent loss functions, such as the piecewise polynomial $l(x) = \frac{1}{2} (x_+)^2$, are not twice differentiable, and thus lack the sufficient conditions based on third-order information. For such functions, in these specific cases, a direct verification approach can be employed to prove the convexity of $H(\rho)$.

\begin{corollary}
Consider the loss function $l(x) = \frac{1}{2} (x_+)^2$. Then the univariate function $H(\rho)$ is convex on $\rho > 0$ under this setting. Furthermore, $H(\rho)$ is differentiable, and the term $\mathcal{K}_H(\rho^k)$ in Line~4 of Algorithm \stmxref{alg:semismooth newton}{3} can be replaced
with $\{H'(\rho^k)\}$, implying that Algorithm \stmxref{alg:semismooth newton}{3} reduces to the classical Newton's method.
The sequence $\{\rho^k\}_{k=0}^\infty$ generated by Algorithm \stmxref{alg:semismooth newton}{3} converges either finitely or $Q$-superlinearly to $\rho^{*}$ satisfying $H(\rho^{*}) = 0$, for any initial point $\rho^0 > 0$.
\end{corollary}

\begin{proof}
We prove the claim by exploiting the explicit formula of $\bs{u}(\rho)$ under $l(x)=\tfrac12(x_+)^2$, and then invoking the convexity and monotonicity of $H$.

From $\bs{G}(\bs{u}, \rho) = \bs{u} - \bs{x} + \frac{\rho}{m} \nabla L(\bs{u}) = \bs{0}$, the $i{\text{th}}$ component is
\begin{equation*}
u_i(\rho) - x_i + \frac{\rho}{m} l'(u_i(\rho)) = u_i(\rho) - x_i + \frac{\rho}{m} (u_i(\rho))_+ = 0.
\end{equation*}
Since $\rho>0$ with the fixed $x_i$, it is easy to show that $u_i(\rho)$ has the explicit form
\begin{equation*}
u_i(\rho) = \left\{
\begin{aligned}
&\frac{mx_i}{m + \rho} & \text{if} \quad x_i > 0,\\
& x_i & \text{if} \quad x_i \le 0.
\end{aligned}
\right.
\end{equation*}
Then $H(\rho)$ has the explicit form
\begin{equation}
\xlabel{eq:rho}
H(\rho)=\sum_{i=1}^m \tfrac12\,(u_i(\rho)_+)^2 - m\lambda
=\sum_{i=1}^m\frac12\left(\frac{m}{m+\rho} (x_i)_+\right)^2 - m\lambda,
\end{equation}
which is differentiable on $(0,\infty)$ with
\begin{equation*}
H'(\rho)=-\sum_{i=1}^m \frac{m^2 (x_i)_+^2}{(m+\rho)^3}.
\end{equation*}
It is easy to observe that $H'(\rho) < 0$, which means that $H(\rho)$ is monotonically decreasing. Moreover, it is obvious that $H(\rho)$ is convex in $\rho$ for $\rho > 0$ from its form in \eqref{eq:rho}. Similar to the proof of  \xref{lem:H well-defined}, it is easy to see that $H$ satisfies the nondegeneracy condition, i.e.,~$H'(\rho^{*})<0$. The rest of the proof follows similarly to that of Theorem 2.1 in~\textcite{aravkin2019level}.
\end{proof}

\section{Numerical Experiments}\xlabel{sec:experiments}
In this section, we conduct experiments on a Linux server equipped with 256~GB RAM and a 96-core AMD EPYC~7402 CPU running at 2.8~GHz.

If the loss function is twice continuously differentiable, the problem~\eqref{p:projection on Z} can be solved using the classical primal--dual interior-point method. However, due to the sparsity of the problem, most of the matrices generated during the iterations are diagonal, which may result in redundant computations when using a standard interior-point solver. To address this, we develop an interior-point method tailored to exploit the sparsity structure of problem~\eqref{p:projection on Z}, as detailed in Section~B.3 of the supplementary material. In the projection experiments, following the experimental setup of~\textcite{luxenberg2025operator}, we compare our algorithm with the implemented interior-point method, the bisection method from Algorithm \stmxref{alg:bisec projection on Z}{4}, the classical solver MOSEK~\parencite{mosek}, and Clarabel~\parencite{Clarabel_2024}. Clarabel is a numerical interior-point solver for convex optimization problems, designed for quadratic objective functions subject to linear and conic constraints. Then we compare our algorithm with the classical solver MOSEK and an alternative splitting-based algorithm, the Splitting Conic Solver (SCS)~\parencite{scs}, to evaluate both the effectiveness and computational efficiency of the proposed method in practical scenarios. Our algorithm is implemented in Python~3.8, with both comparative solvers invoked via CVXPY~\parencite{diamond2016cvxpy}. For each experiment, we set a 1000-second time limit for each method.

We consider two data generation schemes.

\textbf{Synthetic Data.} Following the setup in~\textcite{hu2018utility} and \textcite{hegde2024online}, we let $r_i$ denote the return of the $i$th underlying asset. The expected returns $\E[r_i]$ are uniformly distributed between 0.05 and 0.50, with monotonically increasing values for $i = 1,\ldots,n$. The standard deviation is defined as $\text{std}[r_i] = \E[r_i] + 0.05$, $i=1,\ldots,n$. The covariance between any pair of distinct returns $(r_i,r_j)$ with $i\ne j$ is set to $0.35\,\text{std}[r_i]\text{std}[r_j]$.

\textbf{Real-world Data.} The first dataset contains 458 risky assets from the S\&P 500 index, with data spanning from August 31, 2009, to August 20, 2013. The second dataset comprises 300 constituent stocks of the CSI 300 Index as of January 2010, covering the period from January 4, 2010, to December 31, 2024. The third dataset includes 1000 stocks from the CSI 1000 Index as of January 2015, with data ranging from January 5, 2015, to December 31, 2024. All data sources and detailed preprocessing steps are provided in Section~B.2 of the supplementary material.

For the loss function in UBSR, we consider two types, namely, the exponential function, defined as $l(x) = \exp(\beta x)$ with $\beta > 0$, and the piecewise polynomial function, defined as $l(x) = \eta^{-1}(x_+)^{\eta}$ with $\eta > 1$. Further experimental details are provided in Section~B.1 of the supplementary material.

\subsection{{\?{UBSR}} Projection}

We generate random UBSR projection problems using the vector $\bs{x} \in \mathbb{R}^m$, with entries following a standard normal distribution. We consider both the exponential loss function and the piecewise polynomial loss function, along with various combinations of $(\beta, \lambda)$ or $(\eta, \lambda)$. For each parameter pair, we conduct five independent randomized trials (i.e.,~five distinct random seeds are used to generate $\bs{x}$) to ensure statistical reliability.

As shown in \xref{fig: projection exp}, ``DirSSN\_proj'' and ``SepSSN\_proj'' perform comparably and consistently outperform existing solvers across all settings. Moreover, both methods generally surpass ``IPM\_proj'' and ``Bi\_proj'' in almost all cases, especially when the loss function is a piecewise polynomial function (where for $\eta = 2$, the second derivative of $l$ in the interior-point method is replaced by a subgradient of $\nabla l$). Across all tested dimensions, their runtime remains below 10~s, and in high-dimensional scenarios, they exhibit speedups of more than three orders of magnitude compared with existing solvers. It is worth noting that Clarabel is also a solver based on the interior-point method.
We can observe that the specially implemented interior-point method tailored to the problem structure achieves higher efficiency,
yet it still performs somewhat worse than the proposed method.
These results demonstrate the robustness and effectiveness of the proposed methods.

However, it is worth noting that both ``DirSSN\_proj'' and ``IPM\_proj'' require different initialization strategies under varying loss functions and parameter configurations to achieve the reported convergence performance, highlighting the sensitivity of ``DirSSN\_proj'' and ``IPM\_proj'' to initialization. In practice, we carefully tuned the parameters for ``IPM\_proj'', which require the adjustment of multiple parameters. We emphasize that, for fixed parameters, it is challenging to achieve good results across all settings. In contrast, all the results for ``SepSSN\_proj'' shown in the figures were obtained using a fixed initialization of $\rho_0 = 1$, indicating that ``SepSSN\_proj'' is insensitive to the choice of initialization. In practice, ``IPM\_proj'' performs slightly worse compared to the proposed methods. Clearly, ``Bi\_proj'' also performs poorly. Thus, ``SepSSN\_proj'' consistently delivers robust high performance in practice. Moreover, the loss functions employed in these numerical experiments ensure the global convergence of the G-semismooth Newton method, which typically exhibits faster practical convergence. Therefore, we adopt Algorithm \stmxref{alg:semismooth newton}{3} for problem \eqref{p:projection on Z} in the subsequent experiments.

\begin{figure}[t]
\caption{\xlabel{fig: projection exp}Comparison of the problem dimensions and computational time for solving projection problem~\eqref{p:projection on Z} under different settings. Solid lines and filled markers correspond to $\lambda=0.1$, while dashed lines and empty markers correspond to $\lambda=0.2$.}
\includegraphics[width=0.85\linewidth]{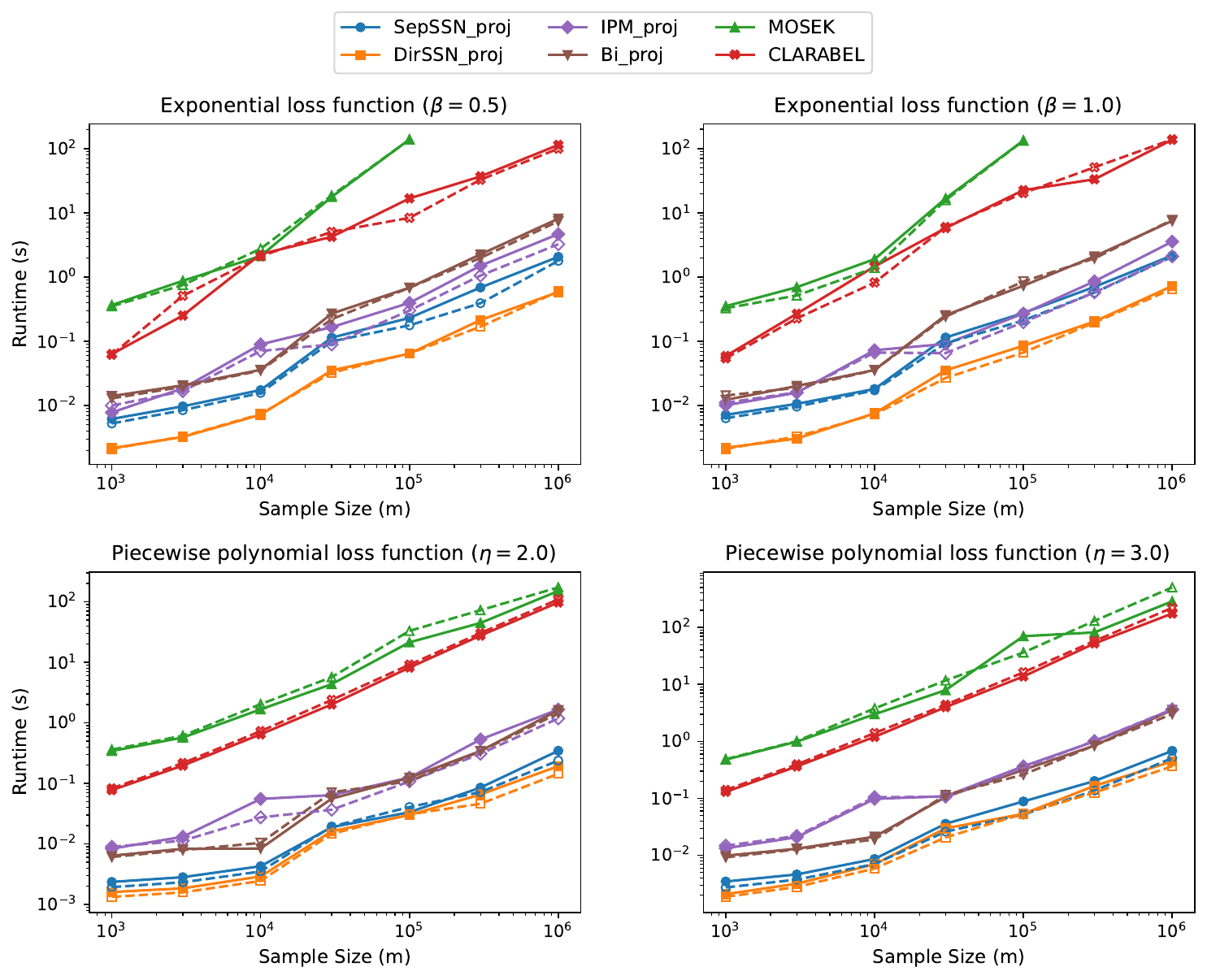}
\end{figure}

\subsection{Comparison with Stochastic Gradient Methods under a \?{Gauss}ian Setting}

Stochastic gradient (SG) methods are often attractive for large-scale risk minimization problems, but they become less convenient when the UBSR portfolio optimization problem involves additional constraints, such as the expected-return constraint considered here, since each stochastic update then requires a projection onto a more complicated feasible set. In contrast, the proposed ADMM method handles these equality and inequality constraints naturally through the SAA reformulation and solves the resulting constrained problem efficiently. To further assess solution quality and enable a direct comparison with SG methods, we consider a Gaussian setting for which the UBSR admits a closed-form expression. Specifically, we use the exponential loss $l(x)=\exp(\beta x)$ with $\beta>0$ and assume that $X\sim \mathcal{N}(\bs{\mu}, \Sigma)$. Then, by \textcite[Example 3.1]{giesecke2008measuring},
\begin{equation*}
\mathrm{SR}_\lambda(X^\top \bs{w})
=
-\bs{\mu}^\top \bs{w}+\frac{\beta}{2}\bs{w}^\top\Sigma \bs{w}-\frac{\log\lambda}{\beta}.
\end{equation*}

Substituting this expression into problem~\eqref{p:original optimization} yields the deterministic convex quadratic program
\begin{equation*}
\min_{\bs{w}\in\mathcal{W}}\quad
\frac{\beta(1-\alpha)}{2}\bs{w}^\top \Sigma \bs{w}
- \bs{\mu}^\top \bs{w}
\quad
\text{s.t.}\quad
\bs{\mu}^\top \bs{w} \ge R_0,
\end{equation*}
up to the constant $-(1-\alpha)\log\lambda/\beta$. Its optimal solution, denoted by $\bs{w}^{*}$, is computed to high precision by a convex QP solver and used as the ground truth. Since $\lambda$ only contributes this constant term, it does not affect the optimal portfolio weight $\bs{w}^{*}$. The weights obtained from the SAA problem solved by ADMM and from the SG methods are similarly insensitive to $\lambda$ in this setting. Therefore, we fix $\lambda=0.1$ and focus on the effects of $(\beta,\alpha)$ and the sample size.

We compare ADMM with the fixed-batch and increasing-batch stochastic gradient (SG) methods of \textcite{gupte2024optimization}. We consider three problem instances with the same number of assets, $n=1000$, and sample sizes $m\in\{5000,20000,50000\}$. For each sample size and random seed, all methods use the same pool of $m$ Gaussian return samples. ADMM constructs the SAA problem using the entire sample pool. Here, the sampling budget counts distinct observations; deterministic passes over the fixed SAA sample do not constitute additional sampling. To impose the same $m$-sample budget on the SG methods, we randomly permute the sample pool and partition it without replacement across 200 iterations. Let $b_k$ denote the number of samples allocated to iteration $k$, where $\sum_{k=1}^{200}b_k=m$, so that each sample is assigned to exactly one iteration. At iteration $k$, all $b_k$ samples are used to estimate the mean return required for the projection step. We then divide these samples into two disjoint subbatches of nearly equal size: one is used to estimate the UBSR value $t_k$, and the other is used to estimate the stochastic gradient based on the resulting estimate of $t_k$.

For fixed-batch SG, we set $b_k=m/200$ for all $k$. For increasing-batch SG, we use the nondecreasing schedule
\fontsize{7pt}{10.5pt}\selectfont
\begin{equation*}
\widetilde b_k=2+\left\lfloor\frac{(m-400)k}{\sum_{j=1}^{200}j}\right\rfloor,\qquad
g=m-\sum_{k=1}^{200}\widetilde b_k,\qquad
b_k=\widetilde b_k+\mathbf{1}_{\{k>200-g\}},\quad k=1,\ldots,200.
\end{equation*}\normalsize
Thus, each iteration receives at least two observations, the remaining budget is allocated in proportion to $1,2,\ldots,200$, and any rounding remainder is assigned to the final $g$ iterations. This preserves monotonicity and ensures $\sum_{k=1}^{200}b_k=m$. Both SG methods use the stepsize $\eta_k=10/k$. The reported runtime for the SG methods excludes the generation of the common sample pool, whereas the runtime for ADMM includes it. We report the CPU time, absolute objective gap, and solution error $\| \bs{w}-\bs{w}^{*}\| $, each averaged over 20 independent random seeds.

\tabref{tab:ubsr_comp} reports the results. Compared with fixed-batch SG, \ubrk increasing-batch SG attains a smaller objective gap in every setting and a smaller $\| \bs{w}-\bs{w}^{*}\| $ in all but one setting, at a modest additional time cost. ADMM attains the smallest objective gap and is faster than both SG variants in every setting. Both error measures decrease as $m$ increases. Overall, ADMM obtains the highest objective accuracy with the lowest computational cost under the same sample budget.

\begin{table}[t]
\caption{\xlabel{tab:ubsr_comp}Comparison of computational performance between ADMM and SG methods under matched total sampling budgets in the Gaussian setting. ``SG (fixed)'' denotes fixed-batch SG, and ``SG (increasing)'' denotes increasing-batch SG.}

\resizebox{\textwidth}{!}{
\begin{tabular}{LLLLLL}
\toprule
Total sample budget $m$ & $(\beta,\alpha)$ & Method & Avg. Time (s) & Avg. Obj. Gap & Avg. $\| \boldsymbol{w}-\boldsymbol{w}^*\| $ \\
\midrule
\multirow{12}*{5000}
& \multirow{3}*{$(0.5,0.0)$} & SG (fixed) & 18.16 &2.27e\tmin02 &2.26e\tmin01 \\
& & SG (increasing) & 18.99 &1.48e\tmin02 &2.13e\tmin01 \\
& & ADMM & 5.36 &2.94e\tmin03 &2.30e\tmin01 \\
\cmidrule{2-6}
& \multirow{3}*{$(0.5,0.5)$} & SG (fixed) & 15.47 &1.72e\tmin02 &2.64e\tmin01 \\
& & SG (increasing) & 17.66 &1.06e\tmin02 &2.46e\tmin01 \\
& & ADMM & 3.83 &3.16e\tmin03 &3.27e\tmin01 \\
\cmidrule{2-6}
& \multirow{3}*{$(1.0,0.0)$} & SG (fixed) & 17.73 &2.22e\tmin02 &1.88e\tmin01 \\
& & SG (increasing) & 18.62 &1.49e\tmin02 &1.74e\tmin01 \\
& & ADMM & 5.75 &2.77e\tmin03 &1.60e\tmin01 \\
\cmidrule{2-6}
& \multirow{3}*{$(1.0,0.5)$} & SG (fixed) & 16.08 &1.68e\tmin02 &2.17e\tmin01 \\
& & SG (increasing) & 18.04 &1.03e\tmin02 &1.98e\tmin01 \\
& & ADMM & 4.78 &2.96e\tmin03 &2.30e\tmin01 \\
\hline
\multirow{12}*{20,000}
& \multirow{3}*{$(0.5,0.0)$} & SG (fixed) & 15.28 &7.14e\tmin03 &1.68e\tmin01 \\
& & SG (increasing) & 17.48 &3.90e\tmin03 &1.55e\tmin01 \\
& & ADMM & 8.27 &9.40e\tmin04 &1.33e\tmin01 \\
\cmidrule{2-6}
& \multirow{3}*{$(0.5,0.5)$} & SG (fixed) & 16.91 &4.28e\tmin03 &1.86e\tmin01 \\
& & SG (increasing) & 18.84 &2.59e\tmin03 &1.83e\tmin01 \\
& & ADMM & 5.70 &1.10e\tmin03 &2.02e\tmin01 \\
\cmidrule{2-6}
& \multirow{3}*{$(1.0,0.0)$} & SG (fixed) & 15.17 &7.83e\tmin03 &1.39e\tmin01 \\
& & SG (increasing) & 17.49 &4.10e\tmin03 &1.20e\tmin01 \\
& & ADMM & 8.61 &7.74e\tmin04 &8.79e\tmin02 \\
\cmidrule{2-6}
& \multirow{3}*{$(1.0,0.5)$} & SG (fixed) & 18.20 &4.30e\tmin03 &1.43e\tmin01 \\
& & SG (increasing) & 19.91 &2.37e\tmin03 &1.32e\tmin01 \\
& & ADMM & 7.89 &9.43e\tmin04 &1.34e\tmin01 \\
\hline
\multirow{12}*{50,000}
& \multirow{3}*{$(0.5,0.0)$} & SG (fixed) & 17.24 &2.54e\tmin03 &1.25e\tmin01 \\
& & SG (increasing) & 19.56 &1.40e\tmin03 &1.08e\tmin01 \\
& & ADMM & 11.37 &3.89e\tmin04 &8.60e\tmin02 \\
\cmidrule{2-6}
& \multirow{3}*{$(0.5,0.5)$} & SG (fixed) & 18.24 &1.52e\tmin03 &1.40e\tmin01 \\
& & SG (increasing) & 18.89 &1.08e\tmin03 &1.43e\tmin01 \\
& & ADMM & 8.71 &4.95e\tmin04 &1.36e\tmin01 \\
\cmidrule{2-6}
& \multirow{3}*{$(1.0,0.0)$} & SG (fixed) & 14.66 &2.84e\tmin03 &1.02e\tmin01 \\
& & SG (increasing) & 17.32 &1.43e\tmin03 &7.92e\tmin02 \\
& & ADMM & 14.29 &3.24e\tmin04 &5.69e\tmin02 \\
\cmidrule{2-6}
& \multirow{3}*{$(1.0,0.5)$} & SG (fixed) & 19.56 &1.48e\tmin03 &1.00e\tmin01 \\
& & SG (increasing) & 20.69 &8.58e\tmin04 &9.26e\tmin02 \\
& & ADMM & 13.23 &3.90e\tmin04 &8.61e\tmin02 \\
\bottomrule
\end{tabular}
}
\end{table}

\subsection{Algorithm Block-wise Behavior Analysis}

To illustrate the computational behavior of the proposed ADMM method, we present a small-scale numerical example with $m=5000$ and $n=1000$. This experiment is designed to provide a more detailed view of the convergence pattern and the computational burden of each ADMM block.

\xref{fig:convergence_curve_response} illustrates the convergence behavior of ADMM in terms of CPU time. Across both the exponential and piecewise polynomial loss functions, the relative objective error decreases rapidly within the first few seconds, indicating that the proposed method is able to obtain high-accuracy solutions efficiently. The convergence curves are also stable across different choices of $(\lambda,\beta)$ or $(\lambda,\eta)$, suggesting that the algorithm is not overly sensitive to different parameters.

\tabref{tab:small_scale_response} further reports the block-wise computational cost of ADMM. The results show that most of the runtime is spent on the $(\bs{w},t)$-update, while the $(\bs{z},s)$-update accounts for a smaller portion of the total time, providing a clear breakdown of the computational burden across the two primal blocks. Moreover, the final primal and dual residuals remain at comparable magnitudes across different loss functions and parameter settings. This indicates that the dynamic residual adjustment mechanism effectively balances primal feasibility and dual feasibility during the iterations.

Overall, this example indicates that the proposed ADMM framework is robust to the choice of loss function and parameter setting. It reaches the prescribed primal and dual residual tolerances within 30--57 iterations in the reported runs, suggesting fast convergence and reliable primal--dual feasibility.

\begin{figure}[t]
\caption{\xlabel{fig:convergence_curve_response}Convergence behavior (relative objective error versus CPU time) for the small-scale numerical example based on a representative single random seed.}
\includegraphics[width=0.8\linewidth]{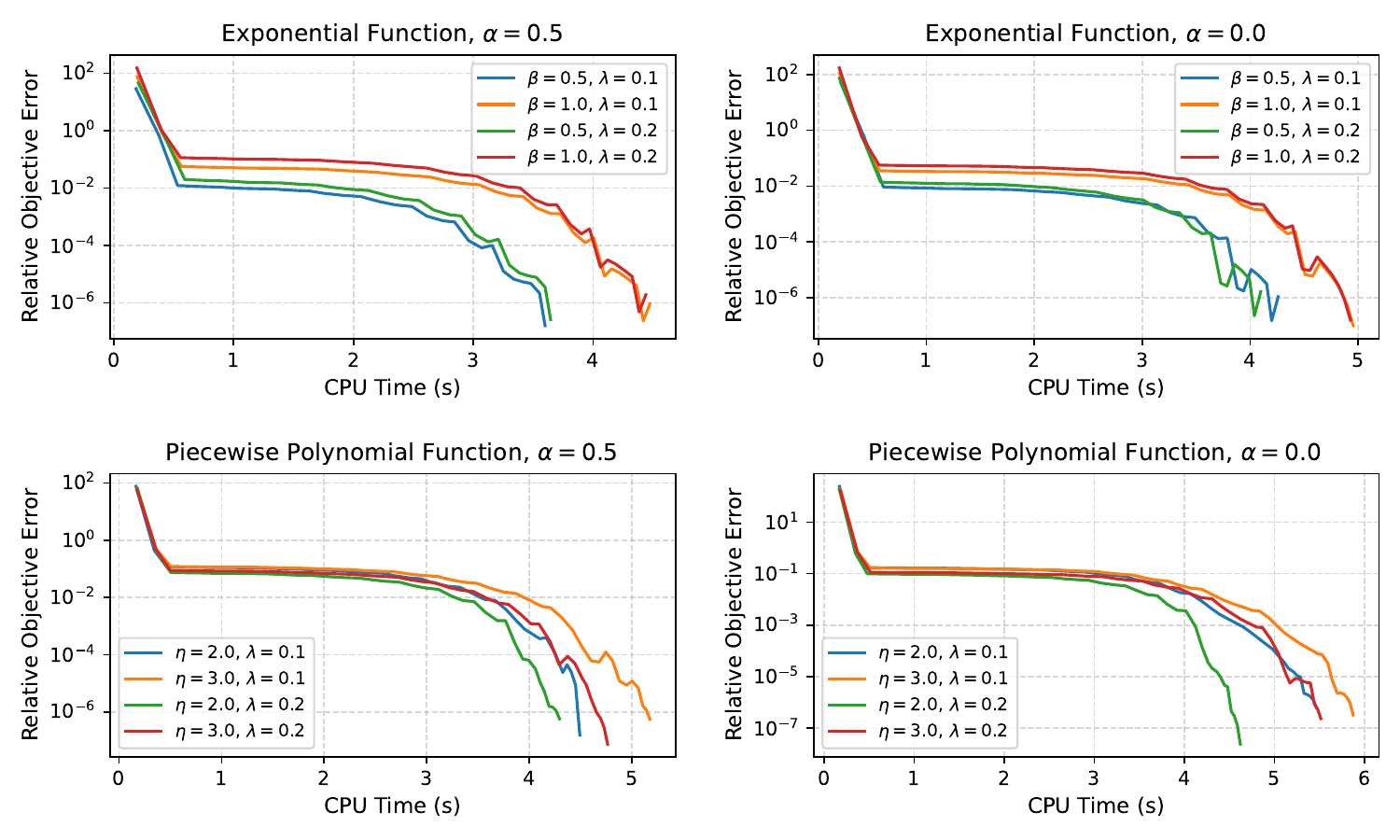}
\end{figure}

\begin{table}[t]
\caption{\xlabel{tab:small_scale_response}Detailed computational breakdown for the small-scale experiment, averaged over 5 random seeds; standard deviations are reported in parentheses.}

\resizebox{\textwidth}{!}{
\begin{tabular}{LLLLLL}
\toprule
{$(\beta/\eta, \lambda)$} & {$\alpha$} & {time/iter} & {$(\boldsymbol{w},t)$-{time}/$(\boldsymbol{z},s)$-{time}} & {{Pri. res}} & {{Dual res}} \\

&  & (mean (std)) & (mean (std)) &  &  \\
\midrule
\multicolumn{6}{C}{ {exponential loss}} \\
\hline
\multirow{2}*{(0.5, 0.1)} & 0.5 & 3.24 (0.44)/30.0 (0.0) & 2.80 (0.34)/0.36 (0.09) &1.88e\tmin04(6.33e\tmin06) &1.87e\tmin05(8.74e\tmin07) \\
& 0.0 & 4.05 (0.13)/34.0 (0.0) & 3.37 (0.13)/0.59 (0.01) &1.48e\tmin04(9.06e\tmin06) &1.68e\tmin05(1.81e\tmin06) \\
\multirow{2}*{(0.5, 0.2)} & 0.5 & 3.53 (0.12)/30.0 (0.0) & 2.99 (0.12)/0.46 (0.00) &1.87e\tmin04(6.32e\tmin06) &1.87e\tmin05(8.75e\tmin07) \\
& 0.0 & 4.10 (0.14)/34.0 (0.0) & 3.41 (0.14)/0.58 (0.00) &1.48e\tmin04(9.06e\tmin06) &1.68e\tmin05(1.81e\tmin06) \\
\multirow{2}*{(1.0, 0.1)} & 0.5 & 4.42 (0.21)/37.0 (0.0) & 3.74 (0.16)/0.58 (0.07) &1.39e\tmin04(3.29e\tmin05) &1.53e\tmin05(6.85e\tmin06) \\
& 0.0 & 5.00 (0.19)/41.2 (0.4) & 4.14 (0.20)/0.74 (0.01) &1.50e\tmin04(4.71e\tmin05) &1.85e\tmin05(8.94e\tmin06) \\
\multirow{2}*{(1.0, 0.2)} & 0.5 & 4.37 (0.21)/37.0 (0.0) & 3.68 (0.22)/0.59 (0.01) &1.39e\tmin04(3.29e\tmin05) &1.53e\tmin05(6.85e\tmin06) \\
& 0.0 & 5.00 (0.21)/41.4 (0.5) & 4.15 (0.20)/0.73 (0.01) &1.22e\tmin04(2.59e\tmin05) &1.38e\tmin05(8.46e\tmin06) \\
\hline
\multicolumn{6}{C}{{piecewise polynomial loss}} \\
\hline
\multirow{2}*{(2.0, 0.1)} & 0.5 & 4.66 (0.09)/42.0 (0.7) & 4.43 (0.08)/0.12 (0.00) &5.79e\tmin05(1.02e\tmin05) &2.38e\tmin05(1.10e\tmin05) \\
& 0.0 & 5.60 (0.14)/56.6 (0.5) & 5.29 (0.15)/0.16 (0.01) &8.21e\tmin05(8.42e\tmin06) &1.15e\tmin05(6.10e\tmin07) \\
\multirow{2}*{(2.0, 0.2)} & 0.5 & 4.19 (0.05)/40.0 (0.0) & 3.98 (0.06)/0.10 (0.01) &7.74e\tmin05(1.44e\tmin05) &9.31e\tmin06(2.06e\tmin06) \\
& 0.0 & 4.72 (0.07)/45.2 (0.4) & 4.45 (0.08)/0.13 (0.01) &9.85e\tmin05(1.56e\tmin05) &1.87e\tmin05(4.15e\tmin06) \\
\multirow{2}*{(3.0, 0.1)} & 0.5 & 5.20 (0.21)/47.4 (0.5) & 4.76 (0.23)/0.31 (0.01) &7.92e\tmin05(1.83e\tmin05) &1.40e\tmin05(5.93e\tmin06) \\
& 0.0 & 6.04 (0.23)/56.8 (0.8) & 5.48 (0.24)/0.40 (0.01) &8.16e\tmin05(1.69e\tmin05) &1.50e\tmin05(1.22e\tmin05) \\
\multirow{2}*{(3.0, 0.2)} & 0.5 & 4.80 (0.22)/43.8 (0.4) & 4.40 (0.22)/0.28 (0.01) &1.05e\tmin04(1.20e\tmin05) &8.40e\tmin06(5.15e\tmin06) \\
& 0.0 & 5.53 (0.22)/51.8 (0.8) & 5.02 (0.22)/0.37 (0.01) &8.89e\tmin05(1.62e\tmin05) &1.23e\tmin05(7.22e\tmin06) \\
\bottomrule
\end{tabular}
}
\end{table}

\subsection{{\?{UBSR}} Optimization}

Based on the SAA technique, we solve problem~\eqref{p:saa optimization}. Accordingly, we conduct experiments on problem~\eqref{p:saa optimization} using both synthetic and real-world datasets. We examine various $(\beta, \lambda)$ or $(\eta, \lambda)$ combinations and perform five independent random replications for each parameter pair to ensure statistical reliability. The threshold $R_0$ is set to the expected future return under the $1/n$ policy.

\textbf{Results.} The results of the synthetic data are shown in \tabref{tab: syn data exp,tab: syn data poly}, while the results for the real-world data are presented in \tabref{tab: real data exp,tab: real data poly}, with all reported values representing the average of five independent runs with different random seeds.
``ADMM'' represents the results of our proposed algorithm, ``MOSEK'' refers to the results obtained from the MOSEK solver, and ``SCS'' denotes the results from the Splitting Conic Solver. ``obj'' indicates the objective value of the SAA problem~\eqref{p:saa optimization}, and ``time'' reports the corresponding solution time in seconds. ``/'' indicates that the constraint violation defined in~\eqref{eq:constraint violation} exceeded 10\tsup{\tmin5} in at least three out of five experiments. Bold numbers highlight the best results. The parameter $\alpha$ represents the $\alpha$ value of the objective function in problem~\eqref{p:saa optimization}. Note that when $\alpha=0$, it corresponds to the classical portfolio optimization problem considered in \textcite{hegde2024online,hu2018utility}. The constraint violation for problem~\eqref{p:saa optimization} is defined as follows:
\begin{align}
\xlabel{eq:constraint violation}
\nf&\max\left\{\vbox to 1.2pc{}\max\{[w_i]_{-},i=1,\ldots,n\},  \svert \bs{1}_n^\top \bs{w} - 1\svert ,  [\bs{\mu}^\top \bs{w} - R_0]_{-},\right.\nonumber\\&\xq \left.\left[\lambda - \frac{1}{m}\sum_{i=1}^m l(-\xi_i^\top\bs{w}-t)\right]_{-}\right\}.
\end{align}

\begin{table}[t]
\caption{\xlabel{tab: syn data exp}Comparison results on synthetic data using exponential loss functions with different dimensions for problem~\eqref{p:saa optimization}.}
\resizebox{\textwidth}{!}{
\begin{tabular}{LLLLLLLLLLLLLL}
\toprule
$m$& $n$ & \multicolumn{6}{L}{$\alpha=0.5$} & \multicolumn{6}{L}{$\alpha=0$} \\

\cmidrule{3-8} \cmidrule{9-14}
& & \multicolumn{2}{L}{ADMM} & \multicolumn{2}{L}{MOSEK} & \multicolumn{2}{L}{SCS}
& \multicolumn{2}{L}{ADMM} & \multicolumn{2}{L}{MOSEK} & \multicolumn{2}{L}{SCS} \\
\cmidrule{3-4} \cmidrule{5-6} \cmidrule{7-8}
\cmidrule{9-10} \cmidrule{11-12} \cmidrule{13-14}
& & obj & time & obj & time & obj & time & obj & time & obj & time & obj & time \\
\midrule
& & \multicolumn{12}{C}{$\beta = 0.5, \lambda = 0.1$}  \\
\hline
5000 & 500 & 1.8186  & \textbf{1.44 } & 1.8186 & 9.39  &1.8186 & 29.62  &4.1374 & \textbf{1.48 } & 4.1374 & 9.45  &4.1374 & 26.28  \\
5000 & 1000 & 1.8132  & \textbf{3.96 } & 1.8132 & 20.43  &1.8132 & 60.01  &4.1318 & \textbf{4.29 } & 4.1318 & 18.64  &4.1318 & 62.61  \\
30000 & 500 & 1.8213  & \textbf{3.71 } & 1.8213 & 74.60  &/ &/&4.1397 & \textbf{5.18 } & 4.1397 & 73.75  &/ &/ \\
30000 & 1000 & 1.8193  & \textbf{6.22 } & 1.8193 & 128.23  &/ &/&4.1370 & \textbf{8.03 } & 4.1370 & 123.12  &/ &/ \\
100000 & 500 & 1.8215  & \textbf{7.59 } & 1.8215 & 321.83  &/ &/&4.1398 & \textbf{9.36 } & 4.1398 & 317.07  &/ &/ \\
100000 & 1000 & 1.8196  & \textbf{11.44 } & 1.8196 & 557.98  &/ &/&4.1371 & \textbf{13.56 } & 4.1371 & 534.05  &/ &/ \\
\hline
& & \multicolumn{12}{C}{$\beta = 0.5, \lambda = 0.2$}  \\
\hline
5000 & 500 & 1.1255  & \textbf{1.21 } & 1.1255 & 9.62  &1.1255 & 31.75  &2.7511 & \textbf{1.51 } & 2.7511 & 9.05  &2.7511 & 27.76  \\
5000 & 1000 & 1.1200  & \textbf{3.84 } & 1.1200 & 18.99  &1.1200 & 60.76  &2.7455 & \textbf{4.27 } & 2.7455 & 19.37  &2.7455 & 56.51  \\
30000 & 500 & 1.1282  & \textbf{3.78 } & 1.1282 & 72.51  &1.1282 & 267.69  &2.7534 & \textbf{5.04 } & 2.7534 & 69.18  &2.7534 & 212.92  \\
30000 & 1000 & 1.1261  & \textbf{6.11 } & 1.1261 & 135.39  &1.1261 & 528.95  &2.7507 & \textbf{14.16 } & 2.7507 & 140.47  &2.7507 & 581.34  \\
100000 & 500 & 1.1284  & \textbf{7.00 } & 1.1284 & 329.06  &/ &/&2.7535 & \textbf{9.21 } & 2.7535 & 320.98  &/ &/ \\
100000 & 1000 & 1.1265  & \textbf{11.11 } & 1.1265 & 519.77  &/ &/&2.7508 & \textbf{13.66 } & 2.7508 & 518.13  &/ &/ \\
\hline
& & \multicolumn{12}{C}{$\beta = 1, \lambda = 0.1$}  \\
\hline
5000 & 500 & 0.6834  & \textbf{1.74 } & 0.6834 & 10.09  &0.6834 & 37.93  &1.8643 & \textbf{1.92 } & 1.8643 & 9.74  &1.8643 & 38.00  \\
5000 & 1000 & 0.6779  & \textbf{5.00 } & 0.6779 & 21.29  &0.6779 & 182.39  &1.8582 & \textbf{5.31 } & 1.8582 & 20.14  &1.8582 & 94.38  \\
30000 & 500 & 0.6858  & \textbf{5.56 } & 0.6858 & 71.88  &0.6853 & 856.31  &1.8664 & \textbf{6.58 } & 1.8664 & 71.43  &/ &/ \\
30000 & 1000 & 0.6831  & \textbf{8.74 } & 0.6831 & 128.52  &/ &/&1.8631 & \textbf{10.60 } & 1.8631 & 128.71  &/ &/ \\
100000 & 500 & 0.6859  & \textbf{10.26 } & 0.6859 & 402.51  &/ &/&1.8665 & \textbf{14.07 } & 1.8665 & 393.37  &/ &/ \\
100000 & 1000 & 0.6833  & \textbf{15.50 } & 0.6833 & 590.06  &/ &/&1.8630 & \textbf{19.45 } & 1.8630 & 580.71  &/ &/ \\
\hline
& & \multicolumn{12}{C}{$\beta = 1, \lambda = 0.2$}  \\
\hline
5000 & 500 & 0.3369  & \textbf{1.65 } & 0.3369 & 9.62  &0.3369 & 43.16  &1.1711 & \textbf{1.88 } & 1.1711 & 9.88  &1.1711 & 37.77  \\
5000 & 1000 & 0.3313  & \textbf{4.86 } & 0.3313 & 19.55  &0.3313 & 183.19  &1.1650 & \textbf{5.56 } & 1.1650 & 19.94  &1.1650 & 97.25  \\
30000 & 500 & 0.3392  & \textbf{5.41 } & 0.3392 & 72.48  &0.3393 & 389.83  &1.1733 & \textbf{6.66 } & 1.1733 & 73.05  &1.1733 & 256.55  \\
30000 & 1000 & 0.3365  & \textbf{8.90 } & 0.3365 & 128.44  &0.3362 & 711.90  &1.1699 & \textbf{10.35 } & 1.1699 & 129.93  &1.1699 & 747.50  \\
100000 & 500 & 0.3393  & \textbf{11.25 } & 0.3393 & 444.18  &/ &/&1.1733 & \textbf{14.40 } & 1.1733 & 415.94  &/ &/ \\
100000 & 1000 & 0.3367  & \textbf{16.73 } & 0.3367 & 621.48  &/ &/&1.1699 & \textbf{20.90 } & 1.1699 & 604.60  &/ &/ \\
\bottomrule
\end{tabular}
}
\end{table}

\begin{table}[t]
\caption{\xlabel{tab: syn data poly}Comparison results on synthetic data using piecewise polynomial loss functions with different dimensions for problem~\eqref{p:saa optimization}.}

\resizebox{\textwidth}{!}{
\begin{tabular}{LLLLLLLLLLLLLL}
\toprule
$m$&$n$ & \multicolumn{6}{L}{$\alpha=0.5$} & \multicolumn{6}{L}{$\alpha=0$} \\

\cmidrule{3-8} \cmidrule{9-14}
&  & \multicolumn{2}{L}{ADMM} & \multicolumn{2}{L}{MOSEK} & \multicolumn{2}{L}{SCS}
& \multicolumn{2}{L}{ADMM} & \multicolumn{2}{L}{MOSEK} & \multicolumn{2}{L}{SCS} \\
\cmidrule{3-4} \cmidrule{5-6} \cmidrule{7-8}
\cmidrule{9-10} \cmidrule{11-12} \cmidrule{13-14}
& & obj & time & obj & time & obj & time & obj & time & obj & time & obj & time \\
\midrule
& & \multicolumn{12}{C}{$\eta = 2, \lambda = 0.1$}  \\
\hline
5000 & 500 & \tmin 0.6505 & \textbf{1.44 } & \tmin 0.6505 & 8.43  & \tmin 0.6505 & 166.42  & \tmin 0.8092 & \textbf{1.45 } & \tmin 0.8092 & 8.80  & \tmin 0.8092 & 405.70  \\
5000 & 1000 & \tmin 0.6567 & \textbf{4.91 } & \tmin 0.6567 & 17.43  & \tmin 0.6567 & 656.04  & \tmin 0.8159 & \textbf{6.21 } & \tmin 0.8159 & 18.04  & \tmin 0.8159 & 780.92  \\
30000 & 500 & \tmin 0.6483 & \textbf{4.83 } & \tmin 0.6483 & 52.80  &/ &/& \tmin 0.8066 & \textbf{5.13 } & \tmin 0.8066 & 53.94  &/ &/ \\
30000 & 1000 & \tmin 0.6519 & \textbf{9.53 } & \tmin 0.6519 & 105.30  &/ &/& \tmin 0.8104 & \textbf{10.21 } & \tmin 0.8104 & 107.53  &/ &/ \\
100000 & 500 & \tmin 0.6483 & \textbf{12.76 } & \tmin 0.6483 & 190.56  &/ &/& \tmin 0.8063 & \textbf{13.35 } & \tmin 0.8063 & 183.88  &/ &/ \\
100000 & 1000 & \tmin 0.6519 & \textbf{19.32 } & \tmin 0.6519 & 365.43  &/ &/& \tmin 0.8105 & \textbf{21.64 } & \tmin 0.8105 & 368.00  &/ &/ \\
\hline
& & \multicolumn{12}{C}{$\eta = 2, \lambda = 0.2$}  \\
\hline
5000 & 500 & \tmin 0.7636 & \textbf{1.23 } & \tmin 0.7636 & 8.16  & \tmin 0.7636 & 189.01  & \tmin 1.0329 & \textbf{1.41 } & \tmin 1.0329 & 8.48  & \tmin 1.0329 & 691.32  \\
5000 & 1000 & \tmin 0.7697 & \textbf{4.50 } & \tmin 0.7697 & 18.38  &/ &/& \tmin 1.0395 & \textbf{4.90 } & \tmin 1.0395 & 18.38  &/ &/ \\
30000 & 500 & \tmin 0.7615 & \textbf{4.21 } & \tmin 0.7615 & 54.00  &/ &/& \tmin 1.0307 & \textbf{4.80 } & \tmin 1.0307 & 51.82  &/ &/ \\
30000 & 1000 & \tmin 0.7648 & \textbf{7.81 } & \tmin 0.7648 & 103.34  &/ &/& \tmin 1.0346 & \textbf{9.33 } & \tmin 1.0346 & 106.12  &/ &/ \\
100000 & 500 & \tmin 0.7616 & \textbf{9.37 } & \tmin 0.7616 & 181.80  &/ &/& \tmin 1.0307 & \textbf{12.00 } & \tmin 1.0307 & 181.12  &/ &/ \\
100000 & 1000 & \tmin 0.7648 & \textbf{15.72 } & \tmin 0.7648 & 378.85  &/ &/& \tmin 1.0347 & \textbf{18.91 } & \tmin 1.0347 & 356.19  &/ &/ \\
\hline
& & \multicolumn{12}{C}{$\eta = 3, \lambda = 0.1$}  \\
\hline
5000 & 500 & \tmin 0.7473 & \textbf{1.59 } & \tmin 0.7473 & 9.95  &/ &/& \tmin 1.0045 & \textbf{1.78 } & \tmin 1.0045 & 9.58  &/ &/ \\
5000 & 1000 & \tmin 0.7536 & \textbf{5.45 } & \tmin 0.7536 & 18.42  &/ &/& \tmin 1.0112 & \textbf{6.72 } & \tmin 1.0112 & 18.98  &/ &/ \\
30000 & 500 & \tmin 0.7449 & \textbf{6.26 } & \tmin 0.7449 & 55.65  &/ &/& \tmin 1.0012 & \textbf{6.70 } & \tmin 1.0012 & 56.21  &/ &/ \\
30000 & 1000 & \tmin 0.7485 & \textbf{10.93 } & \tmin 0.7485 & 108.13  &/ &/& \tmin 1.0048 & \textbf{12.24 } & \tmin 1.0048 & 102.30  &/ &/ \\
100000 & 500 & \tmin 0.7448 & \textbf{16.07 } & \tmin 0.7448 & 202.18  &/ &/& \tmin 1.0006 & \textbf{16.98 } & \tmin 1.0006 & 198.60  &/ &/ \\
100000 & 1000 & \tmin 0.7486 & \textbf{23.21 } & \tmin 0.7486 & 404.06  &/ &/& \tmin 1.0048 & \textbf{26.35 } & \tmin 1.0048 & 377.18  &/ &/ \\
\hline
& & \multicolumn{12}{C}{$\eta = 3, \lambda = 0.2$}  \\
\hline
5000 & 500 & \tmin 0.8502 & \textbf{1.39 } & \tmin 0.8502 & 9.33  & \tmin 0.8503 & 898.46  & \tmin 1.2083 & \textbf{1.61 } & \tmin 1.2083 & 9.75  &/ &/ \\
5000 & 1000 & \tmin 0.8564 & \textbf{5.09 } & \tmin 0.8564 & 17.92  &/ &/& \tmin 1.2150 & \textbf{6.65 } & \tmin 1.2150 & 17.98  &/ &/ \\
30000 & 500 & \tmin 0.8480 & \textbf{5.62 } & \tmin 0.8480 & 57.94  &/ &/& \tmin 1.2055 & \textbf{5.79 } & \tmin 1.2055 & 56.32  &/ &/ \\
30000 & 1000 & \tmin 0.8515 & \textbf{10.13 } & \tmin 0.8515 & 101.99  &/ &/& \tmin 1.2094 & \textbf{11.65 } & \tmin 1.2094 & 102.20  &/ &/ \\
100000 & 500 & \tmin 0.8479 & \textbf{14.05 } & \tmin 0.8479 & 188.38  &/ &/& \tmin 1.2052 & \textbf{16.01 } & \tmin 1.2052 & 187.03  &/ &/ \\
100000 & 1000 & \tmin 0.8515 & \textbf{20.24 } & \tmin 0.8515 & 388.44  &/ &/& \tmin 1.2095 & \textbf{25.58 } & \tmin 1.2095 & 363.32  &/ &/ \\
\bottomrule
\end{tabular}
}
\end{table}

\begin{table}[t]
\caption{\xlabel{tab: real data exp}Comparison results on real-world data using exponential loss functions with different dimensions for problem~\eqref{p:saa optimization}.}

\resizebox{\textwidth}{!}{
\begin{tabular}{LLLLLLLLLLLLLL}
\toprule
$m$& $n$& \multicolumn{6}{L}{$\alpha=0.5$} & \multicolumn{6}{L}{$\alpha=0$} \\

\cmidrule{3-8} \cmidrule{9-14}
&  & \multicolumn{2}{L}{ADMM} & \multicolumn{2}{L}{MOSEK} & \multicolumn{2}{L}{SCS}
& \multicolumn{2}{L}{ADMM} & \multicolumn{2}{L}{MOSEK} & \multicolumn{2}{L}{SCS} \\
\cmidrule{3-4} \cmidrule{5-6} \cmidrule{7-8}
\cmidrule{9-10} \cmidrule{11-12} \cmidrule{13-14}
& & obj & time & obj & time & obj & time
& obj & time & obj & time & obj & time \\
\midrule
& & \multicolumn{12}{C}{$\beta = 0.5, \lambda = 0.1$}  \\
\hline
800 & 458 & 2.2997  & \textbf{0.28 } & 2.2997 & 8.95  &2.2997 & 1.45  &4.6024 & \textbf{0.32 } & 4.6024 & 10.31  &4.6024 & 1.38  \\
1000 & 458 & 2.2999  & \textbf{0.48 } & 2.2999 & 12.49  &2.2999 & 1.75  &4.6026 & \textbf{0.52 } & 4.6026 & 13.06  &4.6026 & 1.93  \\
2000 & 300 & 2.3010  & \textbf{0.24 } & 2.3010 & 11.13  &2.3010 & 2.49  &4.6036 & \textbf{0.35 } & 4.6036 & 11.53  &4.6036 & 1.92  \\
3000 & 300 & 2.3011  & \textbf{0.26 } & 2.3011 & 15.77  &2.3011 & 3.64  &4.6038 & \textbf{0.32 } & 4.6038 & 15.93  &4.6038 & 3.73  \\
1000 & 1000 & 2.2988  & \textbf{0.76 } & 2.2988 & 32.23  &2.2988 & 7.89  &4.6014 & \textbf{0.69 } & 4.6014 & 32.59  &4.6014 & 8.13  \\
2000 & 1000 & 2.2988  & \textbf{0.67 } & 2.2988 & 64.00  &2.2988 & 20.97  &4.6015 & \textbf{0.88 } & 4.6015 & 64.83  &4.6015 & 14.86  \\
\hline
& & \multicolumn{12}{C}{$\beta = 0.5, \lambda = 0.2$}  \\
\hline
800 & 458 & 1.6066  & \textbf{0.32 } & 1.6066 & 9.99  &1.6066 & 1.40  &3.2161 & \textbf{0.43 } & 3.2161 & 9.65  &3.2161 & 1.56  \\
1000 & 458 & 1.6067  & \textbf{0.59 } & 1.6067 & 12.76  &1.6067 & 1.68  &3.2163 & \textbf{0.48 } & 3.2163 & 12.64  &3.2163 & 2.14  \\
2000 & 300 & 1.6079  & \textbf{0.25 } & 1.6079 & 11.41  &1.6079 & 2.86  &3.2174 & \textbf{0.38 } & 3.2174 & 10.78  &3.2174 & 2.04  \\
3000 & 300 & 1.6080  & \textbf{0.25 } & 1.6080 & 16.60  &1.6080 & 3.65  &3.2175 & \textbf{0.34 } & 3.2175 & 16.40  &3.2175 & 2.81  \\
1000 & 1000 & 1.6057  & \textbf{0.68 } & 1.6057 & 32.76  &1.6057 & 8.32  &3.2151 & \textbf{0.76 } & 3.2151 & 32.37  &3.2151 & 7.92  \\
2000 & 1000 & 1.6057  & \textbf{0.86 } & 1.6057 & 64.17  &1.6057 & 21.73  &3.2152 & \textbf{0.75 } & 3.2152 & 64.46  &3.2152 & 15.33  \\
\hline
& & \multicolumn{12}{C}{$\beta = 1, \lambda = 0.1$}  \\
\hline
800 & 458 & 1.1485  & \textbf{0.36 } & 1.1485 & 10.35  &1.1485 & 1.71  &2.3000 & \textbf{1.00 } & 2.3000 & 10.08  &2.3000 & 1.81  \\
1000 & 458 & 1.1487  & \textbf{0.52 } & 1.1487 & 12.59  &1.1487 & 2.48  &2.3002 & \textbf{0.79 } & 2.3002 & 12.78  &2.3002 & 1.67  \\
2000 & 300 & 1.1498  & \textbf{0.43 } & 1.1498 & 12.60  &1.1498 & 2.32  &2.3011 & \textbf{0.43 } & 2.3011 & 12.92  &2.3011 & 2.40  \\
3000 & 300 & 1.1499  & \textbf{0.35 } & 1.1499 & 21.57  &1.1499 & 3.62  &2.3012 & \textbf{0.46 } & 2.3012 & 18.16  &2.3012 & 3.39  \\
1000 & 1000 & 1.1476  & \textbf{0.71 } & 1.1476 & 33.16  &1.1476 & 10.31  &2.2990 & \textbf{1.08 } & 2.2990 & 32.60  &2.2990 & 8.41  \\
2000 & 1000 & 1.1476  & \textbf{0.74 } & 1.1476 & 64.35  &1.1476 & 26.68  &2.2990 & \textbf{0.83 } & 2.2990 & 63.95  &2.2990 & 21.28  \\
\hline
& & \multicolumn{12}{C}{$\beta = 1, \lambda = 0.2$}  \\
\hline
800 & 458 & 0.8020  & \textbf{0.44 } & 0.8020 & 10.52  &0.8020 & 1.97  &1.6069 & \textbf{1.08 } & 1.6069 & 11.15  &1.6069 & 1.39  \\
1000 & 458 & 0.8021  & \textbf{0.55 } & 0.8021 & 13.09  &0.8021 & 2.68  &1.6070 & \textbf{0.85 } & 1.6070 & 13.65  &1.6070 & 2.17  \\
2000 & 300 & 0.8032  & \textbf{0.44 } & 0.8032 & 12.14  &0.8032 & 2.74  &1.6080 & \textbf{0.50 } & 1.6080 & 12.79  &1.6080 & 2.79  \\
3000 & 300 & 0.8033  & \textbf{0.41 } & 0.8033 & 15.88  &0.8033 & 4.64  &1.6081 & \textbf{0.45 } & 1.6081 & 17.22  &1.6081 & 3.63  \\
1000 & 1000 & 0.8010  & \textbf{0.84 } & 0.8010 & 32.67  &0.8010 & 8.95  &1.6058 & \textbf{0.88 } & 1.6058 & 35.09  &1.6058 & 12.45  \\
2000 & 1000 & 0.8010  & \textbf{0.89 } & 0.8010 & 64.28  &0.8010 & 17.64  &1.6058 & \textbf{0.80 } & 1.6058 & 67.97  &1.6058 & 16.04 \\
\bottomrule
\end{tabular}
}
\end{table}

\begin{table}[t]
\caption{\xlabel{tab: real data poly}Comparison results on real-world data using piecewise polynomial loss functions with different dimensions for problem~\eqref{p:saa optimization}.}

\resizebox{\textwidth}{!}{
\begin{tabular}{LLLLLLLLLLLLLL}
\toprule
$m$ &  $n$   & \multicolumn{6}{L}{$\alpha=0.5$} & \multicolumn{6}{L}{$\alpha=0$} \\

\cmidrule{3-8}\cmidrule{9-14}
&  & \multicolumn{2}{L}{ADMM} & \multicolumn{2}{L}{MOSEK} & \multicolumn{2}{L}{SCS}
& \multicolumn{2}{L}{ADMM} & \multicolumn{2}{L}{MOSEK} & \multicolumn{2}{L}{SCS} \\
\cmidrule{3-4}\cmidrule{5-6}\cmidrule{7-8}
\cmidrule{9-10}\cmidrule{11-12}\cmidrule{13-14}
&     & obj & time & obj & time & obj & time & obj & time & obj & time & obj & time \\
\midrule
& & \multicolumn{12}{C}{$\eta = 2, \lambda = 0.1$}  \\
\hline
800 & 458 & \tmin 0.2261 & \textbf{0.69 } & \tmin 0.2261 & 4.72  & \tmin 0.2261 & 3.55  & \tmin 0.4494 & \textbf{2.17 } & \tmin 0.4494 & 5.01  & \tmin 0.4494 & 3.10  \\
1000 & 458 & \tmin 0.2260 & \textbf{0.66 } & \tmin 0.2260 & 5.25  & \tmin 0.2260 & 5.67  & \tmin 0.4493 & \textbf{1.62 } & \tmin 0.4493 & 6.09  & \tmin 0.4493 & 5.41  \\
2000 & 300 & \tmin 0.2250 & \textbf{0.36 } & \tmin 0.2250 & 5.00  & \tmin 0.2250 & 9.70  & \tmin 0.4485 & \textbf{0.52 } & \tmin 0.4485 & 4.89  & \tmin 0.4485 & 7.57  \\
3000 & 300 & \tmin 0.2249 & \textbf{0.39 } & \tmin 0.2249 & 7.58  & \tmin 0.2249 & 20.94  & \tmin 0.4484 & \textbf{0.59 } & \tmin 0.4484 & 7.15  & \tmin 0.4484 & 17.96  \\
1000 & 1000 & \tmin 0.2272 & \textbf{0.80 } & \tmin 0.2272 & 14.89  & \tmin 0.2272 & 11.52  & \tmin 0.4506 & \textbf{1.23 } & \tmin 0.4506 & 15.79  & \tmin 0.4506 & 10.66  \\
2000 & 1000 & \tmin 0.2272 & \textbf{0.81 } & \tmin 0.2272 & 27.51  & \tmin 0.2272 & 43.40  & \tmin 0.4505 & \textbf{0.92 } & \tmin 0.4505 & 29.90  & \tmin 0.4505 & 38.64  \\
\hline
& & \multicolumn{12}{C}{$\eta = 2, \lambda = 0.2$}  \\
\hline
800 & 458 & \tmin 0.3189 & \textbf{0.78 } & \tmin 0.3189 & 5.03  & \tmin 0.3189 & 3.98  & \tmin 0.6348 & \textbf{1.40 } & \tmin 0.6348 & 5.09  & \tmin 0.6348 & 3.40  \\
1000 & 458 & \tmin 0.3187 & \textbf{0.52 } & \tmin 0.3187 & 6.39  & \tmin 0.3187 & 9.25  & \tmin 0.6347 & \textbf{1.03 } & \tmin 0.6347 & 6.19  & \tmin 0.6347 & 6.47  \\
2000 & 300 & \tmin 0.3177 & \textbf{0.32 } & \tmin 0.3177 & 5.73  & \tmin 0.3177 & 12.72  & \tmin 0.6338 & \textbf{0.57 } & \tmin 0.6338 & 5.29  & \tmin 0.6338 & 9.98  \\
3000 & 300 & \tmin 0.3176 & \textbf{0.27 } & \tmin 0.3176 & 8.07  & \tmin 0.3176 & 22.00  & \tmin 0.6337 & \textbf{0.64 } & \tmin 0.6337 & 8.06  & \tmin 0.6337 & 19.69  \\
1000 & 1000 & \tmin 0.3199 & \textbf{0.75 } & \tmin 0.3199 & 14.84  & \tmin 0.3199 & 12.00  & \tmin 0.6359 & \textbf{0.89 } & \tmin 0.6359 & 16.37  & \tmin 0.6359 & 11.58  \\
2000 & 1000 & \tmin 0.3199 & \textbf{0.98 } & \tmin 0.3199 & 29.81  & \tmin 0.3199 & 54.89  & \tmin 0.6359 & \textbf{0.88 } & \tmin 0.6359 & 34.30  & \tmin 0.6359 & 49.68  \\
\hline
& & \multicolumn{12}{C}{$\eta = 3, \lambda = 0.1$}  \\
\hline
800 & 458 & \tmin 0.3371 & \textbf{1.21 } & \tmin 0.3371 & 4.89  & \tmin 0.3371 & 5.54  & \tmin 0.6714 & \textbf{1.97 } & \tmin 0.6714 & 4.84  & \tmin 0.6714 & 4.46  \\
1000 & 458 & \tmin 0.3370 & \textbf{0.79 } & \tmin 0.3370 & 6.25  & \tmin 0.3370 & 8.21  & \tmin 0.6713 & \textbf{1.72 } & \tmin 0.6713 & 5.88  & \tmin 0.6713 & 6.24  \\
2000 & 300 & \tmin 0.3361 & \textbf{0.58 } & \tmin 0.3361 & 6.10  & \tmin 0.3361 & 25.04  & \tmin 0.6707 & \textbf{0.71 } & \tmin 0.6707 & 5.46  & \tmin 0.6707 & 18.97  \\
3000 & 300 & \tmin 0.3360 & \textbf{0.57 } & \tmin 0.3360 & 10.41  & \tmin 0.3360 & 48.61  & \tmin 0.6706 & \textbf{0.80 } & \tmin 0.6706 & 7.97  & \tmin 0.6706 & 44.30  \\
1000 & 1000 & \tmin 0.3382 & \textbf{0.92 } & \tmin 0.3382 & 6.00  & \tmin 0.3382 & 27.49  & \tmin 0.6726 & \textbf{0.94 } & \tmin 0.6726 & 6.71  & \tmin 0.6726 & 24.42  \\
2000 & 1000 & \tmin 0.3382 & \textbf{0.95 } & \tmin 0.3382 & 29.48  & \tmin 0.3382 & 165.85  & \tmin 0.6726 & \textbf{1.13 } & \tmin 0.6726 & 30.88  & \tmin 0.6726 & 117.54  \\
\hline
& & \multicolumn{12}{C}{$\eta = 3, \lambda = 0.2$}  \\
\hline
800 & 458 & \tmin 0.4242 & \textbf{0.70 } & \tmin 0.4242 & 4.74  & \tmin 0.4242 & 8.88  & \tmin 0.8456 & \textbf{2.33 } & \tmin 0.8456 & 4.89  & \tmin 0.8456 & 7.06  \\
1000 & 458 & \tmin 0.4241 & \textbf{0.83 } & \tmin 0.4241 & 5.55  & \tmin 0.4241 & 10.70  & \tmin 0.8455 & \textbf{1.70 } & \tmin 0.8455 & 5.71  & \tmin 0.8455 & 8.06  \\
2000 & 300 & \tmin 0.4231 & \textbf{0.44 } & \tmin 0.4231 & 6.25  & \tmin 0.4231 & 39.32  & \tmin 0.8447 & \textbf{0.54 } & \tmin 0.8447 & 5.43  & \tmin 0.8447 & 28.20  \\
3000 & 300 & \tmin 0.4230 & \textbf{0.46 } & \tmin 0.4230 & 8.54  & \tmin 0.4230 & 76.47  & \tmin 0.8446 & \textbf{0.64 } & \tmin 0.8446 & 8.96  & \tmin 0.8446 & 73.34  \\
1000 & 1000 & \tmin 0.4253 & \textbf{1.71 } & \tmin 0.4253 & 5.89  & \tmin 0.4253 & 54.24  & \tmin 0.8467 & \textbf{1.10 } & \tmin 0.8468 & 6.46  & \tmin 0.8468 & 39.70  \\
2000 & 1000 & \tmin 0.4253 & \textbf{0.89 } & \tmin 0.4253 & 31.45  & \tmin 0.4253 & 277.42  & \tmin 0.8467 & \textbf{1.04 } & \tmin 0.8467 & 30.65  & \tmin 0.8467 & 264.04 \\
\bottomrule
\end{tabular}
}

\end{table}

As shown in \tabref{tab: syn data exp,tab: syn data poly}, even when using the smoother exponential loss function, the splitting conic solver SCS fails in high-dimensional settings. This highlights the superiority of the ADMM algorithm we developed for portfolio optimization problems with UBSR as the risk measure. Additionally, in some cases, the objective function value of SCS may be lower, which we attribute to some instances where SCS produces solutions that violate the constraints, i.e.,~SCS provides incorrect solutions. The table clearly indicates that, across all parameter settings, our algorithm consistently achieves the highest computational efficiency, with its advantages becoming even more pronounced as the sample size $m$ increases. Specifically, for the exponential loss function, MOSEK takes more than five times as long as ADMM on small-scale problems and nearly 50 times as long when the sample size reaches 100,000. For the piecewise polynomial loss function, ADMM is nearly 20 times faster than MOSEK in large-scale settings and more than 100 times faster than SCS. This trend emphasizes the increasing computational advantage of our method as the problem size grows, demonstrating the promising effectiveness of our proposed algorithm for large-scale problems.

On the other hand, from \tabref{tab: real data exp,tab: real data poly}, we observe that SCS performs well in low-dimensional settings and typically outperforms MOSEK when UBSR uses the exponential loss function. This trend reverses for the piecewise polynomial loss function. Nevertheless, our ADMM algorithm consistently achieves shorter computational times across all settings, providing solutions in nearly 1~s for almost all test cases. These results highlight the practical efficiency of our method and its robustness in practical applications.

\subsection{Empirical Investment Performance}

Similar to the empirical setting in~\textcite{cui2024decision}, we evaluate the empirical performance of the proposed UBSR portfolio model. The results are reported in Section~C of the supplementary material.

\section{Conclusion}

This paper addresses the computational challenges of utility-based shortfall risk (UBSR) optimization, a convex risk measure sensitive to tail losses. We propose a novel approach using the Alternating Direction Method of Multipliers (ADMM). A key contribution is the development of two semismooth Newton methods for the projection step in ADMM. Theoretically, we show convergence guarantees for the ADMM and its subproblem solvers. Numerical experiments demonstrate the algorithm's practical benefits on asset management problems with significant speedups over existing solvers. Future work includes incorporating market frictions and transaction costs into the model and extending the model to a distributionally robust setting.

\section*{Acknowledgments}
This work is partly supported by the National Key R\&D Program of China under grant 2023YFA1009300, the National Natural Science Foundation of China under grant 12171100, and the Major Program of the National Natural Science Foundation of China under grants 72394360 and 72394364.

\printbibliography

\clearpage

\setcounter{equation}{16}
\setcounter{algocf}{4}
\setcounter{table}{7}
\setcounter{figure}{3}

\appendix

\section{\texorpdfstring{A Semismooth Newton Method for Solving the KKT System~(9)}{A Semismooth Newton Method for Solving the KKT System}}
\label{sec: dir solve KKT}
A natural approach is to solve the KKT system (9) directly using a semismooth Newton method. Recall that we denote the aggregate loss function as $L(\bs{u}) = \sum^m_{i=1} l(u_i)$, and its gradient is given by $\nabla L(\bs{u}) = [l'(u_1), l'(u_2), \cdots, l'(u_m)]^\top$ where $l'(\cdot)$ denotes the derivative of $l(\cdot)$. Define the mapping $\bs{F}: \mb{R}^{m+1} \to \mb{R}^{m+1}$ by
$$
\bs{F} (\bs{y}) :=
\begin{bmatrix}
	\bs{u} - \bs{x} + \frac{\rho}{m} \nabla L(\bs{u}) \\
	\frac{1}{m} L(\bs{u}) - \lambda
\end{bmatrix},
$$
where $\bs{y}=[\bs{u}^\top, \ \rho ]^\top$. Since all convex functions are semismooth~\parencite[Proposition 3]{mifflin1977semismooth}, under Assumptions~1 and~3, it is easy to see that $\bs{F}$ is semismooth with respect to the following Clarke generalized Jacobian:
\begin{equation}
	\label{eq: JF y}
	\begin{aligned}
		\mathcal{J}_{\bs{F}} (\bs{y}) =  \left\{
		\begin{bmatrix}
			\mf{I} + \frac{\rho}{m} \Lambda & \frac{1}{m} \nabla L(\bs{u}) \\
			\frac{1}{m} \nabla L(\bs{u})^\top & 0
		\end{bmatrix}
		\;\middle|\; \Lambda \in \partial \nabla L(\bs{u})
		\right\},
	\end{aligned}
\end{equation}
where $\partial \nabla L(\bs{u})$ denotes the Clarke subdifferential of $\nabla L(\bs{u})$ and is explicitly given by
\[
\partial \nabla L(\bs{u}) = \left\{ \mathrm{diag}(\eta_1, \eta_2, \ldots, \eta_m) \,\middle|\, \eta_i \in \partial l'(u_i) \text{ for all } i \right\}.
\]

The following lemma characterizes the elements of~\eqref{eq: JF y} and will be used repeatedly.

\begin{lemma}
	\label{lem: Jfy nonsingular}
	Suppose that Assumptions~1 and 3 hold. Then, for any $\Lambda \in \partial\nabla L(\bs{u})$ and $\rho>0$, $\Lambda$ is a diagonal matrix with nonnegative diagonal elements and $\mf{I} + \frac{\rho}{m}\Lambda$ is nonsingular. Moreover, every element $J \in \mathcal{J}_{\bs{F}}(\bs{y})$ is nonsingular if and only if $\nabla L(\bs{u}) \neq \bs{0}$.
\end{lemma}
\begin{proof}
	Since $l$ is convex and differentiable, $l'$ is a nondecreasing function. By Assumption~3, $l'$ is locally Lipschitz continuous because it is semismooth. By \textcite[Proposition 2.3]{jiang1995local}, every $\Lambda \in \partial\nabla L(\bs{u})$ is a diagonal matrix with nonnegative diagonal elements, and hence $\mf{I} + \frac{\rho}{m}\Lambda$ is nonsingular for $\rho>0$. For any element $J \in \mathcal{J}_{\bs{F}}(\bs{y})$, the Schur complement criterion~\parencite[A.5.5]{boyd2004convex} shows that $J$ is nonsingular if and only if its Schur complement $-\frac{1}{m^2}\nabla L(\bs{u})^\top(\mf{I} + \frac{\rho}{m}\Lambda)^{-1}\nabla L(\bs{u})$ is nonzero, which is equivalent to $\nabla L(\bs{u}) \neq \bs{0}$.
\end{proof}

At each iteration, the main computational task of the G-semismooth Newton method is to compute the Newton direction $\bs{d}^{k}$ from $ J^{k}\bs{d}^{k}= - \bs{F}(\bs{y}^{k})$. We next show that the linear system is computationally inexpensive, as it can be solved in $O(m)$. Let $\bs{d}^{k} = [(\bs{d}_1^{k})^\top, d_2^k]^\top$ where $\bs{d}_1^{k} \in \mb{R}^m$ and $d_2^k \in \mb{R}$.
When $\nabla L(\bs{u}^{k}) \neq \textbf{0}$, Lemma~\ref{lem: Jfy nonsingular} shows that the Newton direction can be computed by
\[
\begin{aligned}
	\bs{d}_1^k =  \frac{1}{m}\left(\bs{c}^{k}- d_2^k(\bs{I}+\frac{\rho^{k}}{m}\Lambda^{k})^{-1}\nabla L(\bs{u}^{k})\right), \quad
	d_2^k = \frac{m\left( L(\bs{u}^{k}) - m\lambda\right) + \nabla L(\bs{u}^{k})^\top\bs{c}^{k}}{\nabla L(\bs{u}^{k})^\top(\bs{I}+\frac{\rho^{k}}{m}\Lambda^{k})^{-1}\nabla L(\bs{u}^{k})},
\end{aligned}
\]
where $\bs{c}^{k}=-(\bs{I}+\frac{\rho^{k}}{m}\Lambda^{k})^{-1}\left(m\bs{u}^{k} - m\bs{x} + \rho^{k}\nabla L(\bs{u}^{k})\right)$.
Since $\Lambda^{k}$ is a diagonal matrix, the computational cost of the Newton direction is $\mathcal{O}(m)$. The complete algorithm is presented in Algorithm~\ref{alg:hd semismooth newton}.

In general, a high-dimensional semismooth Newton method with a fixed step size lacks global convergence guarantees when initiated from arbitrary starting points. Therefore, in Lines 7 and 8 of Algorithm~\ref{alg:hd semismooth newton}, we incorporate a line search strategy based on Armijo's rule~\parencite{armijo1966minimization} to enhance robustness. Because Lemma~\ref{lem: Jfy nonsingular} assumes $\rho^k>0$, the step-size search also keeps the multiplier positive. On the other hand, by Lemma~\ref{lem: Jfy nonsingular}, when $\nabla L(\bs{u}^k) = \bs{0} $, the generalized Jacobian $J^k$ is singular, and thus the Newton equation $J^k\bs{d}^k=-\bs{F}(\bs{y}^k)$ may not be solvable.
Recall that $a$ is defined in~(12).
Recall that Assumptions~1 and~3 state that $l$ is nondecreasing, differentiable, and convex, so $l' \geq 0$ and $l'$ is nondecreasing. Thus, under Assumptions~1 and~3, $a$ is well-defined.
Hence we obtain that $\nabla L(\bs{u}^{k})\ne \textbf{0}$ if and only if $\max_i u^k_i > a$.
Therefore, an additional step-size search is provided in Lines 5 and 6 of Algorithm~\ref{alg:hd semismooth newton} to prevent the subsequent iterate from entering this singular regime. However, this step may sometimes slow down the overall algorithm.

\begin{algorithm}[ht!]
	\caption{High-Dimensional Semismooth Newton Method for Problem~(7)}
	\label{alg:hd semismooth newton}
	\DontPrintSemicolon
	\textbf{Input:} Initial point $\bs{y}^{0}=[{(\bs{u}^0)}^{\top}, \ \rho^0 ]^\top$ where $\nabla L(\bs{u}^{0})\ne \textbf{0}$ and $\rho^0>0$, parameters $\sigma \in (0,1)$, $\beta \in (0,1)$.\;
	\For{$k = 0,1,\ldots$}{
		Compute Newton direction $\bs{d}^k = [(\bs{d}_1^{k})^\top, d_2^k]^\top$ such that $J^k \bs{d}^k = -\bs{F}(\bs{y}^k)$ where $J^k \in \mathcal{J}_{\bs{F}}(\bs{y}^k)$ defined in \eqref{eq: JF y}.\;
		Initialize step size $\alpha \gets 1$.\;
		\While{$\rho^k+\alpha d_2^k\leq0$ or $\max_i\{u^k_i+\alpha(d_1^k)_i\} \leq a$}{
			Update $\alpha \gets \beta \alpha$.\;
		}
		\While{$\|\bs{F}(\bs{y}^k + \alpha \bs{d}^k)\| > (1 - \sigma \alpha)\|\bs{F}(\bs{y}^k)\|$}{
			Update $\alpha \gets \beta \alpha$.\;
		}
		Update $\bs{y}^{k+1} \gets \bs{y}^k + \alpha \bs{d}^k$.\;
	}
\end{algorithm}

Next, we prove that there exists a neighborhood $N(\bs{y}^*)$ of the solution point $\bs{y}^*=[\bs{u}^{*\top}, \rho^*]^\top$ in which the Newton sequence is well-defined. That is, $J$ is nonsingular for every $J \in \mathcal{J}_{\bs{F}}(\bs{y})$ and every $\bs{y}\in N(\bs{y}^*)$. This property is of critical importance for establishing the local convergence of the semismooth Newton method.

\begin{proposition}
	\label{pro:nonsingular hd_ssn}
	Suppose that Assumptions~1 and~3 hold and $\bs{x} \notin \mathcal{Z}$ in the problem~(7). Then for any element $J \in \mathcal{J}_{\bs{F}}(\bs{y}^*)$ where $\bs{y}^*$ solves $\bs{F}(\bs{y})=0$, $J$ is nonsingular. Furthermore, there is a neighborhood $N(\bs{y}^*)$ of $\bs{y}^*$ and a constant $C$ such that for any $\bs{y} \in N(\bs{y}^*)$ and any $J \in \mathcal{J}_{\bs{F}}(\bs{y})$, $J$ is nonsingular and $\|J^{-1}\|\leq C$.
\end{proposition}
\begin{proof}
	By Lemma~\ref{lem: Jfy nonsingular}, every $J$ is nonsingular if and only if $\nabla L(\bs{u}^*)\neq \bs{0}$.
	We prove $\nabla L(\bs{u}^*)\neq \bs{0}$ by contradiction. Suppose that $\nabla L(\bs{u}^*)=\bs{0}$. According to the KKT conditions~(9), this leads to $u^*_i = x_i$ for all $i = 1, \dots, m$, and thus $\bs{x} = \bs{u}^* \in \mathcal{Z}$, contradicting the assumption that $\bs{x} \notin \mathcal{Z}$. Therefore, every $J \in \mathcal{J}_{\bs{F}}(\bs{y}^*)$ must be nonsingular. The second statement then follows directly from \textcite[Proposition 3.1]{Qi1993semismooth}.
\end{proof}

This proposition ensures that Algorithm~\ref{alg:hd semismooth newton} is locally well defined.

\begin{theorem}
	\label{thm:convergence_hd_ssn}
	Suppose $\bs{y}^*$ is the unique solution to $\bs{F}(\bs{y}) = 0$. Under Assumptions~1 and~3, there exists $\delta > 0$ such that if $\bs{y}^0$ lies in the neighborhood $\mathcal{N}(\bs{y}^*, \delta)$ of $\bs{y}^*$, the sequence $\{\bs{y}^k\}^\infty_{k=0}$ generated by Algorithm~\ref{alg:hd semismooth newton} converges $Q$-superlinearly to $\bs{y}^*$ and the step size $\alpha$ is always $1$. Moreover, if $\bs{F}(\bs{y})$ is strongly semismooth, we have $\bs{y}^k \to \bs{y}^*$ quadratically.
\end{theorem}

\begin{proof}
	By Proposition~\ref{pro:nonsingular hd_ssn}, the Newton iteration is well-defined locally. The remainder of the proof follows from \textcite[Theorem 3.2]{Qi1993semismooth}.
\end{proof}

However, Algorithm~\ref{alg:hd semismooth newton} still lacks a global convergence guarantee. This motivates the implicit-function approach for solving a univariate equation developed in Section~4 of the paper.
\section{Details for Numerical Experiments}

\subsection{Experimental Setup}
\label{ap: exp setting}

For solver configurations, MOSEK and Clarabel are run with their default settings. For SCS, we set the parameter \texttt{eps} to $10^{-9}$ while keeping all other parameters at their default values. This adjustment ensures that SCS achieves a solution precision comparable to our method and MOSEK.

For the real-world datasets, the sample size $m$ denotes the number of trading days randomly selected from each dataset.

For clarity and convenience, our algorithm occasionally reuses certain symbols (e.g., the letter $ t $) without causing ambiguity. To facilitate reproducibility, we provide the implementation details for Algorithm~1 below. The $\bs{w}$-subproblem is solved using Greedy FISTA~\parencite{liang2022FasterF} to ensure faster convergence. The vector $ \bs{w}_0 $ is initialized with all entries equal to $ 1/n $. The slack variable $s$ is set to 0. The variable $ \bs{z}_0 $ is initialized as $ -\bs{R}\bs{w}_0 - t\bs{1}_m $ where $t$ is set to $0$. The dual variable $ \bs{\nu}_0 $ is set to the zero vector. We update $\sigma$ adaptively. When the ratio of the primal residual to the dual residual exceeds 10, we scale $\sigma$ by a factor of $\tau$ in the appropriate direction, where $\tau$ is set to 2.7 in experiments with real-world data and 1.7 in experiments with synthetic data. The initial value of $\sigma$ is set to $10^{-5}$ in experiments with real-world data and $10^{-6}$ in experiments with synthetic data. The maximum number of iterations is set to $1000$. All runs use the stopping rule
\[
\mathrm{Pri.\ res}_k<\epsilon_{\mathrm{pri}}^k, \qquad \mathrm{Dual\ res}_k<\epsilon_{\mathrm{dual}}^k,
\]
where the primal and dual residuals are
\[
\mathrm{Pri.\ res}_k
=
\left(
\left\|
\bs{z}^k+\bs{R}\bs{w}^k+t^k\mf{1}_m
\right\|^2
+
\left(
\bs{\mu}^{\top}\bs{w}^k-s^k-R_0
\right)^2
\right)^{1/2},
\]
\[
\mathrm{Dual\ res}_k
=
\sigma_k
\left(
\left\|
\bs{R}^{\top}\!\left(\bs{z}^k-\bs{z}^{k-1}\right)
-\bs{\mu}\left(s^k-s^{k-1}\right)
\right\|^2
+
\left|
\mf{1}_m^{\top}\!\left(\bs{z}^k-\bs{z}^{k-1}\right)
\right|^2
\right)^{1/2}.
\]
The corresponding stopping tolerances are
\[
\epsilon_{\mathrm{pri}}^k
=
\sqrt{m+1}\,\epsilon_{\mathrm{abs}}
+\epsilon_{\mathrm{rel}}
\max\!\left\{
\left\|
\begin{bmatrix}
\bs{R}\bs{w}^k+t^k\mf{1}_m\\
\bs{\mu}^{\top}\bs{w}^k
\end{bmatrix}
\right\|,
\left\|
\begin{bmatrix}
\bs{z}^k\\
s^k
\end{bmatrix}
\right\|,
|R_0|
\right\},
\]
\[
\epsilon_{\mathrm{dual}}^k
=
\sqrt{n+1}\,\epsilon_{\mathrm{abs}}
+\epsilon_{\mathrm{rel}}
\left\|
\begin{bmatrix}
\bs{R}^{\top}\bs{\nu}_1^k+\bs{\mu}\nu_2^k\\
\mf{1}_m^{\top}\bs{\nu}_1^k
\end{bmatrix}
\right\|,
\qquad
\epsilon_{\mathrm{abs}}=\epsilon_{\mathrm{rel}}=10^{-6}.
\]
For Section~5.2 of the paper, we set $\epsilon_{\mathrm{abs}}=\epsilon_{\mathrm{rel}}=10^{-4}$.

\subsection{Details for Real-world Data}
\label{ap: exp data}
We describe below the real-world datasets and data processing procedures used in this paper:

\begin{itemize}
	\item The first dataset contains daily returns for 458 constituent stocks of the S\&P 500 Index, spanning the period from August 31, 2009, to August 20, 2013. The original data were obtained from the Center for Research in Security Prices (CRSP) through the Wharton Research Data Services\footnote{\href{https://wrds-www.wharton.upenn.edu}{https://wrds-www.wharton.upenn.edu}} (WRDS) platform.

	\item The second and third datasets were obtained from the China Stock Market \& Accounting Research\footnote{\href{https://data.csmar.com}{https://data.csmar.com}} (CSMAR) database. They consist of daily trading data for A-share-listed companies in the Shanghai and Shenzhen stock markets. The raw data contain missing values resulting from trading suspensions, pre-listing periods, or other market closures. Additionally, a small number of outliers are present, likely caused by extreme volatility following the reopening of stocks after prolonged suspensions. Our data preprocessing involves two steps: first, we replace these outliers with null values; then, we impute all missing values using the corresponding stock's mean return. The resulting datasets are used in our experimental analysis.
\end{itemize}

\subsection{Details for Primal-Dual Interior-Point Method}
\label{sec: Primal-Dual IPM}
As mentioned earlier, when the loss function is twice continuously differentiable,
we can design a primal-dual interior-point method~\parencite[Section 11.7]{boyd2004convex} that is specifically tailored to
exploit the sparse structure of problem~(7).
To apply the primal-dual framework, we construct the Lagrangian function $\mathcal{L}(\bs{u},y)$,
where $y$ is the dual variable associated with the inequality constraint:
\[
\mathcal{L}(\bs{u},y)
= \frac{1}{2}\|\bs{u} - \bs{x}\|^2
+ y (\sum_{i=1}^m l(u_i) - m\lambda ).
\]

The optimal solution $(\bs{u}^*, y^*)$ must satisfy the
Karush--Kuhn--Tucker (KKT) conditions~\parencite[Theorem 3.78]{beck2017first}:
\[
\begin{aligned}
	& \nabla_{\bs{u}} \mathcal{L} = 0
	\;\;\Rightarrow\;\; u_i - x_i + y l'(u_i)  = 0, \quad \forall i=1,\dots,m, \\
	& y(\sum_{i=1}^{m} l(u_i)- m\lambda)=0,
	\quad y \geq 0, \quad \sum_{i=1}^{m} l(u_i)- m\lambda \leq 0.
\end{aligned}
\]

At each interior-point iteration, we solve the nonlinear system
\begin{equation}
	\label{eq: IMP F}
	\mathbf{F}_{t}(\bs{u}, y) :=
	\begin{pmatrix}
		u_1 - x_1 + y l'(u_1) \\
		\vdots \\
		u_m - x_m + y l'(u_m) \\
		-y(\sum_{i=1}^{m} l(u_i)- m\lambda)- 1/t
	\end{pmatrix} = \bs{0},
\end{equation}
where $t = \mu / \bigl(-y(\sum_{i=1}^{m} l(u_i)- m\lambda)\bigr)$ and $\mu > 1$ is the centering parameter.
This system can be solved using Newton’s method.
The Jacobian of $\mathbf{F}_t$ is given by
\begin{equation}
	\label{eq: IMP J}
	J(\bs{u}, y) =
	\begin{pmatrix}
		\operatorname{diag}\!\big(1 + y l''(\bs u)\big) & \nabla L(\bs{u}) \\
		-y\nabla L(\bs{u})^\top & m\lambda-\sum_{i=1}^{m} l(u_i)
	\end{pmatrix},
\end{equation}
where $\operatorname{diag}\!\big(\bs 1 + y l''(\bs u)\big)$ denotes a diagonal matrix with diagonal entries $1 + y l''(u_i)$ for $i=1,\dots,m$.
After solving the linear system
\[
J(\bs{u}^{(k)}, y^{(k)})
\begin{pmatrix}
	{\Delta \bs{u}} \\[2pt]
	\Delta y
\end{pmatrix}
= -\mathbf{F}_{t}(\bs{u}^{(k)}, y^{(k)}),
\]
the iterates are updated as
\[
\bs{u}^{(k+1)} = \bs{u}^{(k)} + \alpha \Delta \bs{u},
\quad y^{(k+1)} = y^{(k)} + \alpha \Delta y,
\]
where $\alpha$ is the step size from the line search.

We summarize the procedure in Algorithm~\ref{alg:primal-dual-ipm}.
For convenience, we omit the subscript $t$ and denote
\[
\mathbf{F}_{\bs{u}} =
\begin{pmatrix}
	u_1 - x_1 + y l'(u_1) \\
	\vdots \\
	u_m - x_m + y l'(u_m)
\end{pmatrix},
\quad
\mathbf{F}_y = -y(\sum_{i=1}^{m} l(u_i)- m\lambda)- 1/t.
\]

By exploiting the block structure of $J(\bs{u},y)$,
after straightforward algebraic manipulation we obtain
\[
\begin{aligned}
	\Delta y &=
	\frac{- \mathbf{F}_y  - y {\nabla L(\bs{u})}^\top D^{-1} \mathbf{F}_{\bs{u}}}
	{y{\nabla L(\bs{u})}^\top D^{-1}\nabla L(\bs{u}) + s}, \\[6pt]
	\Delta \bs{u} &= - D^{-1} \left(\mathbf{F}_{\bs{u}}
	+ \Delta y\nabla L(\bs{u})\right),
\end{aligned}
\]
where $D = \operatorname{diag}\!\big(1 + y l''(u_i)\big)$ and $s=m\lambda - \sum_{i=1}^{m} l(u_i)$.
Thus, only the inversion of the diagonal matrix $D$ is required,
along with a few vector multiplications and scalar operations.
Since these operations have computational complexity $\mathcal{O}(m)$,
the overall cost is substantially lower than directly inverting $J$. Following \textcite[Section 11.7]{boyd2004convex}, we set $\gamma$ to $0.5$, $\nu$ to $0.05$, and $y^{(0)}$ to $10$. For $\mu$, we found that using a fixed value performs poorly in practice. Therefore, we adaptively adjust its value based on the dimensionality of the problem.

\begin{algorithm}[ht!]
	\caption{Primal-Dual Interior-Point Method}
	\label{alg:primal-dual-ipm}
	\DontPrintSemicolon
	\textbf{Initialization:}\; Choose an initial point $(\bs{u}^{(0)},y^{(0)})$ such that $\sum_{i=1}^{m} l(u_i^{(0)}) < m\lambda$ and $y^{(0)} > 0$, and set$~\mu > 1,~\epsilon > 0,~\gamma \in (0,1),~\nu \in(0,1)$. \;
	\For{$k = 0,1,2,\ldots$}{
		Set \( t  = \mu / \bigl(-y^{(k)}(\sum_{i=1}^{m} l(u_i^{(k)})- m\lambda)\bigr) \). \;
		Solve the linear system
		\[
		J(\bs{u}^{(k)}, y^{(k)}) \begin{pmatrix} \Delta \bs{u} \\ \Delta y \end{pmatrix}
		= -\mathbf{F}_{t}(\bs{u}^{(k)}, y^{(k)})
		\]
		to obtain the search direction $(\Delta \bs{u}, \Delta y)$ where $-\mathbf{F}_{t}(\bs{u}^{(k)},y^{(k)})$ is defined in \eqref{eq: IMP F} and $J(\bs{u}^{(k)}, y^{(k)})$ is defined in \eqref{eq: IMP J}. \;
		Initialize $\alpha \gets 1$. \;
		\While{$y^{(k)} + \alpha \Delta y \leq 0$ or $\sum_{i=1}^{m} l(u_i^{(k)}+\alpha\Delta u_i)- m\lambda \geq 0$}{
			Update $\alpha \gets \gamma \alpha$.\;
		}
		\While{$\|\mathbf{F}_t(\bs{u}^{(k)}+\alpha\Delta \bs{u}, y^{(k)} + \alpha \Delta y)\|_2 > (1-\nu\alpha)\|\mathbf{F}_t(\bs{u}^{(k)}, y^{(k)})\|_2$}{
			Update $\alpha \gets \gamma \alpha$.\;
		}
		Update $\bs{u}^{(k+1)} = \bs{u}^{(k)} + \alpha \Delta \bs{u}$ and $y^{(k+1)} = y^{(k)} + \alpha \Delta y$. \;
		\If{$\|\mathbf{F}_t(\bs{u}^{(k+1)}, y^{(k+1)})\|_2 < \epsilon$ and $-y^{(k+1)}(\sum_{i=1}^{m} l(u_i^{(k+1)})- m\lambda) < \epsilon$}{Return $\bs{u}^{(k+1)}$.}
	}
\end{algorithm}

\subsection{Details for UBSR Constraint Problem}
\label{app:max util}

Consider another classical portfolio optimization problem with the UBSR serving as a constraint. This leads to the following convex optimization problem:
$$
\max_{\bs{w} \in \mb{R}^n} \quad \mb{E}[u(X^\top\bs{w})] \qquad
\text{s.t.} \quad  \bs{w} \in \mathcal{W}, ~\text{SR}_\lambda(X^\top\bs{w}) \leq b,
$$
where $u(\cdot)$ is a concave and monotonically increasing utility function. The utility function $u(\cdot)$ is generally distinct from the convex loss function $l(\cdot)$. Using similar SAA techniques, we obtain the following SAA approximation:
\begin{equation}
	\label{p: saa utility optimization problem}
	\min_{\bs{w} \in \mb{R}^n} \quad -U(\bs{R}\bs{w}) = -\sum^m_{i=1}u(\xi_i^\top\bs{w})  \qquad
	\text{s.t.} \quad  \bs{w} \in \mathcal{W},~\frac{1}{m}\sum^m_{i=1} l(-\xi_i^\top \bs{w} - b) \leq \lambda.
\end{equation}

By introducing an auxiliary variable $\bs{z} = -\bs{R} \bs{w} - b \bs{1}_m$, the corresponding augmented Lagrangian function for problem~\eqref{p: saa utility optimization problem} is given by
$$
L_{\sigma}(\bs{w}, \bs{z}; \bs{\nu}) =
- U(\bs{R} \bs{w})
+ \bs{\nu}^\top (\bs{R} \bs{w} + \bs{z} + b \bs{1}_m)
+ \frac{\sigma}{2} \left\| \bs{R} \bs{w} + \bs{z} + b \bs{1}_m \right\|^2,
\quad \bs{w} \in \mathcal{W}, \ \bs{z} \in \mathcal{Z}.
$$
Then we can still apply ADMM to solve problem \eqref{p: saa utility optimization problem}. In particular, the $\bs w$-subproblem still minimizes a smooth objective over $\mathcal{W}$, which can be solved by the accelerated projected gradient method \parencite{nesterov2018lectures}; the $\bs z$-subproblem is a projection onto the set $\mathcal Z$, which can be efficiently solved by methods described in Section~4 of the paper.

\section{Empirical Investment Performance}

Following the empirical setting of~\textcite{cui2024decision}, we evaluate the empirical performance of the proposed UBSR portfolio model~(2) on the Fama--French 48 industry portfolios (FF48), which can be downloaded from Kenneth R. French’s website.
\url{http://mba.tuck.dartmouth.edu/pages/faculty/ken.french/data_library.html}.

The sample covers daily returns from December 14, 2016, to December 1, 2021. The out-of-sample evaluation period runs from December 12, 2017, to December 1, 2021, yielding 1000 daily out-of-sample portfolio realizations. The hyperparameter $\alpha$ is set to 0.3, representing risk aversion. We use a rolling-window procedure with a fixed estimation window of the most recent 250 trading days. For each trading day in the out-of-sample period, we construct scenario returns from the previous 250 trading days. The minimum expected return $R_0$ is set to the mean return across all 48 industry portfolios over the rolling window. All optimization problems are solved using Algorithm~1, with the $\bs{z}$-subproblem solved by Algorithm~3, as implemented in our codebase.

We test two families of loss functions, denoted \texttt{exp} (exponential-type loss) and \texttt{poly} (polynomial-type loss), and vary the tuning parameters $\beta$ and $\lambda$. The explicit forms of the loss functions are listed below:
\[
l(x) = \exp(\beta x), \beta > 0, \quad  l(x) = \tfrac{1}{\eta}(x_+)^\eta, \eta \geq 2.
\]
For each configuration we report the out-of-sample average daily return, the Sharpe ratio, the daily return volatility, and the maximum drawdown observed over the out-of-sample period. Results are shown in Table~\ref{tab:ubsr_empirical}. We also include an equally weighted portfolio as a benchmark to evaluate the performance of our UBSR-based models.

\begin{table}[htbp]
	\centering
	\caption{Out-of-sample performance of UBSR-based portfolios.}
	\label{tab:ubsr_empirical}
	\begin{tabular}{lccrrrr}
		\toprule
		Loss function & $\beta/\eta$ & $\lambda$ & Mean return & Sharpe ratio & Volatility & Max drawdown \\
		\midrule
		exp   & 10 & 1    & 0.0005 & 0.0301 & 0.0162 & $-0.3431$ \\
		exp   & 20 & 1    & 0.0004 & 0.0317 & 0.0135 & $-0.3046$ \\
		poly  & 10 & 0.01 & 0.0005 & 0.0313 & 0.0162 & $-0.3402$ \\
		poly  & 20 & 0.01 & 0.0005 & 0.0357 & 0.0138 & $-0.2963$ \\
		\midrule
		Equal-weighted & N/A & N/A & 0.0006 & 0.0422 & 0.0139 & $-0.3832$ \\
		\bottomrule
	\end{tabular}
\end{table}

From Table~\ref{tab:ubsr_empirical}, the UBSR configurations produce average daily returns on the order of $5\times 10^{-4}$, with Sharpe ratios close to $0.03$ and daily volatilities between approximately $1.35\%$ and $1.62\%$. The equal-weighted benchmark attains the highest average daily return and Sharpe ratio among the reported strategies but also records a larger maximum drawdown than several UBSR variants. Increasing $\beta$ or $\eta$ from 10 to 20 while keeping $\lambda$ fixed tends to reduce realized volatility and improve the Sharpe ratio in our experiments. These observations suggest that (i) the UBSR formulations considered deliver conservative average returns relative to the equal-weighted benchmark, (ii) tuning $\beta$ or $\eta$ can improve risk-adjusted performance by lowering volatility, and (iii) certain UBSR parameterizations can attain materially smaller worst-case drawdowns.

\begin{figure}[htbp]
	\centering
	\includegraphics[width=0.9\textwidth]{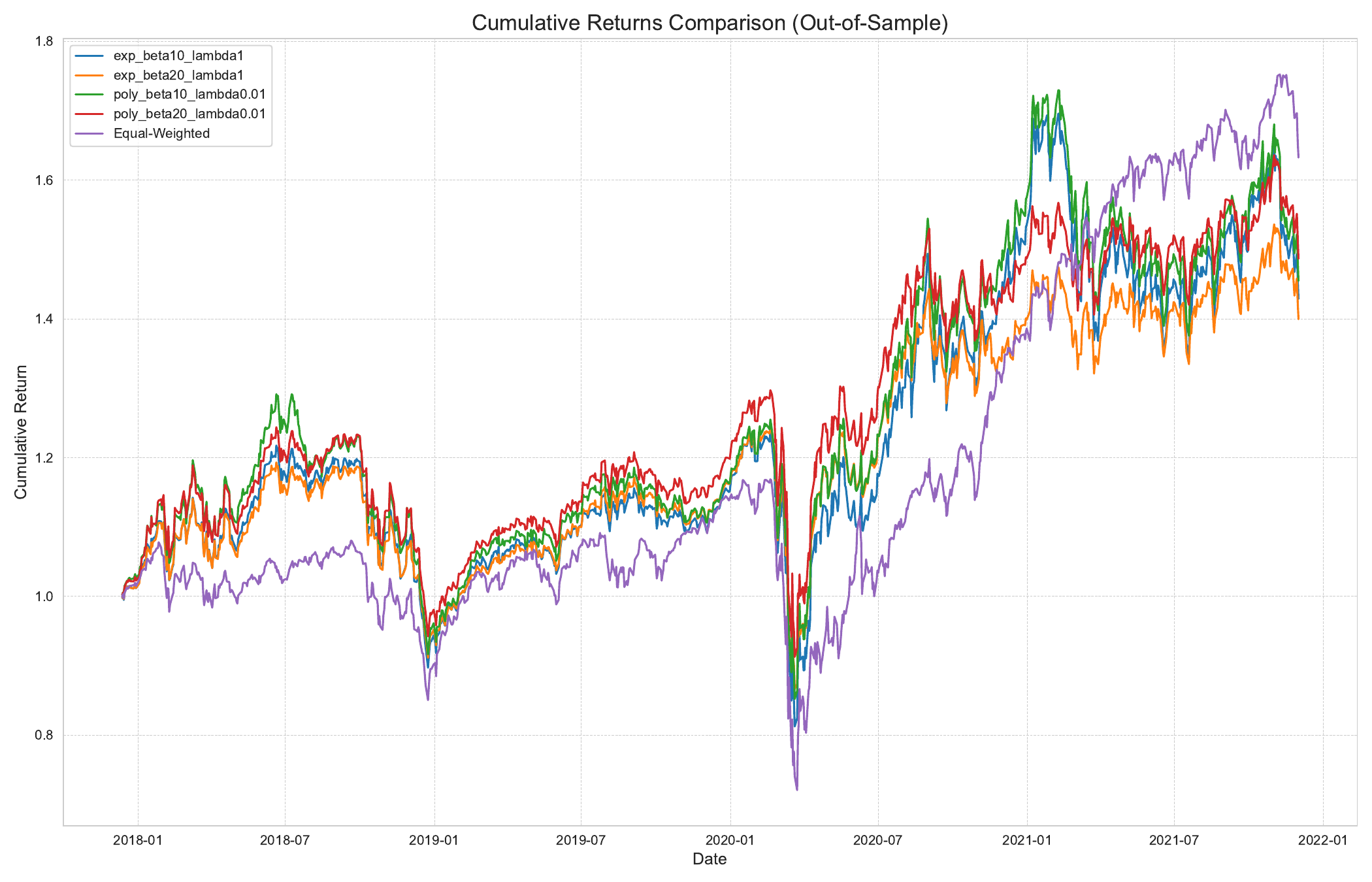}
	\caption{Out-of-sample comparison of cumulative returns.}
	\label{fig:cumulative_returns}
\end{figure}

Figure~\ref{fig:cumulative_returns} presents cumulative returns for the five portfolios over the out-of-sample period.

Before the market downturn in early 2020, several UBSR models (notably \texttt{exp\_beta10\_lambda1} and \texttt{poly\_beta10\_lambda0.01}) outperformed the equal-weighted portfolio, indicating their ability to capture gains in stable markets. During the COVID-19 shock, UBSR models experienced smaller drawdowns than the equal-weighted allocation, suggesting that UBSR-based allocations tend to avoid the most volatile assets. After the decline, the equal-weighted portfolio rebounded more strongly and ultimately produced the highest cumulative return by the end of the sample. This strong rebound contributed substantially to the equal-weighted portfolio's higher terminal return. Nonetheless, several UBSR variants remained competitive by balancing risk control and return generation.

\end{document}